\documentclass[a4paper,12pt,frenchb]{article}

\usepackage{amsmath,amsbsy,amsfonts,amssymb,amsthm}
\usepackage[french]{babel}
\newtheorem*{lem}{Lemme}
\newtheorem*{prop}{Proposition}

\begin{document}

 \title{Fonctions dont les int\'egrales orbitales et celles de leurs transform\'ees de Fourier sont \`a support topologiquement nilpotent  }
\author{J.-L. Waldspurger}
\date{10 octobre 2019}
 \maketitle
 
 {\bf Abstract} Let $F$ be a $p$-adic field and let $G$ be a connected reductive group defined over $F$. We assume $p$ is big. Denote $\mathfrak{g}$ the Lie algebra of $G$. To each vertex $s$ of the reduced Bruhat-Tits' building of $G$, we associate as usual a reductive Lie algebra $\mathfrak{g}_{s}$ defined over the residual field ${\mathbb F}_{q}$. We normalize suitably a Fourier-transform
  $f\mapsto \hat{f}$ on $C_{c}^{\infty}(\mathfrak{g}(F))$.  We study the subspace of functions $f\in C_{c}^{\infty}(\mathfrak{g}(F))$ such that the orbital integrals  of $f$ and of $\hat{f}$  are $0$ for each element of $ \mathfrak{g}(F)$ which  is not topologically nilpotent. This space is related to the characteristic functions of the character-sheaves on the spaces $\mathfrak{g}_{s}({\mathbb F}_{q})$, for each vertex $s$, which are cuspidal and with nilpotent support. We prove that our subspace behave well under endoscopy.
 
 \bigskip 

{\bf Introduction}

\bigskip

Soient $F$ une extension finie d'un corps ${\mathbb Q}_{p}$ et $G$ un groupe r\'eductif connexe d\'efini sur $F$. On impose que $p$ est grand relativement \`a $G$. Introduisons l'immeuble de Bruhat-Tits du groupe adjoint $G_{AD}$. Il est d\'ecompos\'e en facettes et on note $S(G)$ l'ensemble des sommets. A $s\in S(G)$, on associe un sous-groupe parahorique $K_{s}^0\subset G(F)$ et son plus grand sous-groupe distingu\'e pro-$p$-unipotent $K_{s}^+$. D'apr\`es Bruhat et Tits, il existe un groupe r\'eductif connexe $G_{s}$ d\'efini sur le corps r\'esiduel ${\mathbb F}_{q}$ tel que $K_{s}^0/K_{s}^+=G_{s}({\mathbb F}_{q})$. On note par des lettres gothiques les alg\`ebres de Lie de nos diff\'erents groupes. Dans $\mathfrak{g}(F)$, on a de fa\c{c}on similaire une sous-$\mathfrak{o}_{F}$-alg\`ebre $\mathfrak{k}_{s}$ (o\`u $\mathfrak{o}_{F}$ est l'anneau des entiers de $F$) et une sous-$\mathfrak{o}_{F}$-alg\`ebre $\mathfrak{k}_{s}^+$, de sorte que $\mathfrak{k}_{s}/\mathfrak{k}_{s}^+=\mathfrak{g}_{s}({\mathbb F}_{q})$. Toute fonction d\'efinie sur $\mathfrak{g}_{s}({\mathbb F}_{q})$ se rel\`eve en une fonction sur $\mathfrak{k}_{s}$ et s'\'etend par $0$ hors de $\mathfrak{k}_{s}$ en une fonction d\'efinie sur $\mathfrak{g}(F)$. Lusztig a introduit la notion de faisceau-caract\`ere sur $\mathfrak{g}_{s}$. Notons $FC(\mathfrak{g}_{s}({\mathbb F}_{q}))$ l'espace engendr\'e par les fonctions caract\'eristiques des faisceaux-caract\`eres d\'efinis sur $\mathfrak{g}_{s}$ qui sont cuspidaux, \`a support nilpotent et invariants par l'action de Frobenius. Par le proc\'ed\'e ci-dessus, on consid\`ere  $FC(\mathfrak{g}_{s}({\mathbb F}_{q}))$ comme un sous-espace de $C_{c}^{\infty}(\mathfrak{g}(F))$. On note $fc(\mathfrak{g}(F))$ la somme de ces sous-espaces sur tous les sommets $s\in S(G)$. Introduisons l'espace $I(\mathfrak{g}(F))$ qui est le quotient de $C_{c}^{\infty}(\mathfrak{g}(F))$ par le sous-espace des fonctions dont toutes les int\'egrales orbitales sont nulles. L'espace qui nous int\'eresse est l'image $FC(\mathfrak{g}(F))$ de $fc(\mathfrak{g}(F))$ dans $I(\mathfrak{g}(F))$. Pourquoi? Parce que  cet espace joue un r\^ole crucial dans deux travaux ant\'erieurs:  l'un sur les int\'egrales orbitales nilpotentes stables pour les groupes classiques non ramifi\'es, l'autre, avec C. Moeglin, sur les paquets stables de repr\'esentations de r\'eduction unipotente des groupes $SO(2n+1)$. Il nous semble que, pour g\'en\'eraliser ces travaux \`a des groupes quelconques, il est utile de conna\^{\i}tre un peu mieux les propri\'et\'es de cet espace $FC(\mathfrak{g}(F))$, en particulier son comportement par endoscopie. 

Notre premier r\'esultat caract\'erise cet espace $FC(\mathfrak{g}(F))$. On sait que l'on peut d\'efinir une transformation de Fourier $f\mapsto \hat{f}$ dans $C_{c}^{\infty}(\mathfrak{g}(F))$ de sorte que, en notant ${\bf 1}_{k_{s}}$ et ${\bf 1}_{k_{s}^+}$ les fonctions caract\'eristiques de $k_{s}$ et $k_{s}^+$, la fonction $\hat{{\bf 1}}_{k_{s}}$ soit proportionnelle \`a ${\bf 1}_{k_{s}^+}$ pour tout $s\in S(G)$. Cette transformation de Fourier se descend \`a l'espace $I(\mathfrak{g}(F))$. Notons $I_{cusp}(\mathfrak{g}(F))$ le sous-espace des fonctions $f\in I(\mathfrak{g}(F))$ telles que les int\'egrales orbitales $I^G(X,f)$ soient nulles pour tout \'el\'ement $X\in \mathfrak{g}(F)$ qui est semi-simple r\'egulier et non elliptique.

\begin{prop}{L'espace $FC(\mathfrak{g}(F))$ co\"{\i}ncide avec celui des $f\in I(\mathfrak{g}(F))$ telles que les int\'egrales orbitales $I^G(X,f)$ soient nulles pour tout \'el\'ement $X\in \mathfrak{g}(F)$ qui est semi-simple r\'egulier et non topologiquement nilpotent et telles que les int\'egrales orbitales $I^G(X,\hat{f})$ v\'erifient la m\^eme propri\'et\'e. De plus, $FC(\mathfrak{g}(F))$ est contenu dans l'espace $I_{cusp}(\mathfrak{g}(F))$.}\end{prop}

Cette proposition est l'analogue du r\'esultat de Lusztig (\cite{L}) sur les groupes finis. 

La th\'eorie de l'endoscopie est assez simple quand on la restreint \`a l'espace $I_{cusp}(\mathfrak{g}(F))$. Notons $Endo_{ell}(G)$ l'ensemble des classes d'\'equivalence de donn\'ees endoscopiques elliptiques de $G$ (cf. par exemple \cite{MW} chapitre I; une donn\'ee endoscopique ${\bf G}'$ d\'etermine un groupe endoscopique $G'$ mais la notion de donn\'ee endoscopique est plus fine que celle de groupe endoscopique).  On a une d\'ecomposition
$$I_{cusp}(\mathfrak{g}(F))=\oplus_{{\bf G}'}I_{cusp}(\mathfrak{g}(F),{\bf G}'),$$
o\`u $I_{cusp}(\mathfrak{g}(F),{\bf G}')$ est le sous-espace des $f\in I_{cusp}(\mathfrak{g}(F))$ dont les transferts \`a toute donn\'ee ${\bf G}''\not={\bf G}'$ sont nuls. Supposons $G$ quasi-d\'eploy\'e. Il y a une donn\'ee endoscopique elliptique principale ${\bf G}$, pour laquelle le groupe endoscopique est $G$ lui-m\^eme. On note $I_{cusp}^{st}(\mathfrak{g}(F))=I_{cusp}(\mathfrak{g}(F),{\bf G})$. Revenons au cas g\'en\'eral. Pour ${\bf G}'\in Endo_{ell}(G)$, le transfert se restreint en un isomorphisme
$$transfert^{{\bf G}'}: I_{cusp}(\mathfrak{g}(F),{\bf G}')\to I_{cusp}^{st}(\mathfrak{g}'(F))^{Aut({\bf G}')},$$
o\`u, \`a droite, il s'agit du sous-espace des invariants par l'action du groupe fini $Aut({\bf G}')$ des automorphismes de ${\bf G}'$. 

Pour tout ${\bf G}'\in Endo_{ell}(G)$, posons $FC(\mathfrak{g}(F),{\bf G}')=FC(\mathfrak{g}(F))\cap I_{cusp}(\mathfrak{g}(F),{\bf G}')$. Dans le cas o\`u $G$ est quasi-d\'eploy\'e, on pose simplement $FC^{st}(\mathfrak{g}(F))=FC(\mathfrak{g}(F),{\bf G})$. Notre second r\'esultat exprime que les espaces $FC(\mathfrak{g}(F))$ se comportent "bien" par endoscopie. Pr\'ecis\'ement

 \begin{prop}{(i) On a l'\'egalit\'e
 $$FC(\mathfrak{g}(F))=\oplus_{{\bf G}'\in Endo_{ell}(G)}FC(\mathfrak{g}(F),{\bf G}').$$
 
 (ii) Supposons $G$ quasi-d\'eploy\'e. L'espace $FC^{st}(\mathfrak{g}(F))$ est invariant par tout automorphisme de $G$.  
  
 (iii) Pour ${\bf G}'\in Endo_{ell}(G)$, le transfert se restreint en un isomorphisme de $FC(\mathfrak{g}(F),{\bf G}')$ sur $FC^{st}(\mathfrak{g}'(F))^{Aut({\bf G}')}$.}\end{prop}
 
 Dans le paragraphe \ref{produitscalaire}, nous d\'efinirons une application antilin\'eaire injective $D^G:FC(\mathfrak{g}(F))\to I(\mathfrak{g}(F))^*$,  ce dernier espace \'etant celui des distributions sur $\mathfrak{g}(F)$ invariantes par conjugaison par $G(F)$. En renfor\c{c}ant l'hypoth\`ese sur $p$, nous montrerons au paragraphe \ref{traduction} que les espaces images $D^G(FC(\mathfrak{g}(F)))$ v\'erifient des propri\'et\'es analogues \`a celles de la proposition ci-dessus.

 \section{Notations}\label{notations}
 
 Soit $G$ un groupe agissant  sur un ensemble $E$. On note $E^G$ l'ensemble des points fixes. Pour $e\in E$, on note $Z_{G}(e)$ le fixateur de $e$ dans $G$. Pour un sous-ensemble $E'\subset E$, on note $Norm_{G}(E')$ le stabilisateur de $E'$ dans $G$.  
 
 Soient $k$ un corps commutatif parfait. On note $\bar{k}$ une cl\^oture alg\'ebrique de $k$. Pour toute extension $k'$ de $k$ contenue dans $\bar{k}$, on note $\Gamma_{k'/k}$ le groupe de Galois de $k'/k$. On pose simplement  $\Gamma_{k}=\Gamma_{\bar{k}/k}$.  

 Si $E$ est un espace vectoriel sur $k$, on note $E^*$ son dual. Pour tout groupe alg\'ebrique $H$ d\'efini sur $k$, on note $H^0$ sa composante neutre. Soit $T$ un tore d\'efini sur $k$. On note $X_{*}(T)$ le groupe des sous-groupes \`a un param\`etre de $T$, $X^*(T)$ le groupe des caract\`eres de $T$ et on pose ${\cal T}=X_{*}(T)\otimes_{{\mathbb Z}}{\mathbb R}$. Remarquons que $X^*(T)$ s'identifie \`a un ensemble de formes lin\'eaires sur ${\cal T}$. 

Soit $G$ un groupe r\'eductif connexe d\'efini sur $k$.  On note $\mathfrak{g}$ l'alg\`ebre de Lie de $G$. On identifie $G$ et $\mathfrak{g}$ \`a leurs ensembles de points $G(\bar{k})$ et $\mathfrak{g}(\bar{k})$ d\'efinis sur $\bar{k}$. 
On appelle conjugaison l'action adjointe de $G$ sur $\mathfrak{g}$ et on la note selon le cas
 $ (g,X)\mapsto gXg^{-1} $
ou
$  (g,X)\mapsto Ad(g)(X) $. Pour une fonction $f$ sur $\mathfrak{g}(k)$ et pour $g\in G(k)$, on note $^gf$ la fonction $X\mapsto f(g^{-1}Xg)$. 

On note $Z(G)$ le centre de $G$,   $G_{AD}$ le groupe adjoint de $G$ et $G_{SC}$ le rev\^etement simplement connexe du groupe d\'eriv\'e de $G$.  Notons $h(G)$ le plus grand nombre de Coxeter des composantes irr\'eductibles de $G_{AD}$. Si $k$ est de caract\'eristique $p>0$, on suppose $p>h(G)$. 

 On a la d\'ecomposition canonique
$$\mathfrak{g}=\mathfrak{z}(G)\oplus \mathfrak{g}_{SC}.$$
On note  $\mathfrak{g}_{nil}$ le sous-ensemble des \'el\'ements nilpotents de $\mathfrak{g}$. Tous les sous-groupes alg\'ebriques de $G$ seront implicitement suppos\'es d\'efinis sur $k$, sauf mention explicite du contraire. Pour un sous-groupe parabolique $P$ de $G$, on note $U_{P}$ son radical unipotent. On appelle groupe de Levi de $G$ toute composante de Levi d'un sous-groupe parabolique de $G$ (tous deux d\'efinis sur $k$).

 Il arrivera que l'on d\'efinisse un objet relatif au groupe $G$ sans faire figurer la lettre $G$ dans la notation. Quand on aura besoin du m\^eme objet relatif \`a un autre groupe $H$, on ajoutera la lettre $H$ dans la notation.

\section{Groupes sur un corps fini}\label{corpsfinis}
Soit $q$ une puissance enti\`ere d'un nombre premier $p$. On note ${\mathbb F}_{q}$ le  corps fini \`a $q$ \'el\'ements. On note $Fr$ l'\'el\'ement de Frobenius qui engendre $Gal(\bar{{\mathbb F}}_{q}/{\mathbb F}_{q})$. 

Soit  
$G$ un groupe r\'eductif connexe d\'efini sur ${\mathbb F}_{q}$. Comme dans le paragraphe pr\'ec\'edent, on suppose $p>h(G)$.

On note $C(\mathfrak{g}({\mathbb F}_{q}))$ l'espace des fonctions sur $\mathfrak{g}({\mathbb F}_{q})$, \`a valeurs complexes.  Il est muni du produit hermitien non d\'eg\'en\'er\'e
$$(f',f)=\vert G({\mathbb F}_{q})\vert ^{-1}\sum_{X\in \mathfrak{g}({\mathbb F}_{q}}\bar{f}'(X)f(X).$$
On note  $C_{inv}(\mathfrak{g}({\mathbb F}_{q}))$ le sous-espace des fonctions qui sont invariantes par conjugaison par $G({\mathbb F}_{q})$. 

(1) {\bf  Remarque.} La notation est impr\'ecise car $G_{SC}$ et $G_{AD}$ ont la m\^eme alg\`ebre de Lie mais  l'invariance  par conjugaison par $G_{SC}({\mathbb F}_{q})$ est plus faible que l'invariance par $G_{AD}({\mathbb F}_{q})$. Quand il y aura une ambigu\"{\i}t\'e sur le groupe en question, on notera plut\^ot $C_{inv}^G(\mathfrak{g}({\mathbb F}_{q}))$ l'espace ci-dessus. La m\^eme remarque s'applique \`a diff\'erents espaces de fonctions associ\'es  dans la suite \`a des alg\`ebres de Lie, par exemple l'espace $I(\mathfrak{g}(F))$ du paragraphe \ref{groupespadiques} ou l'espace $FC(\mathfrak{g}(F))$ du paragraphe  \ref{FC}.

\bigskip

On note $C_{cusp}(\mathfrak{g}({\mathbb F}_{q}))$ le sous-espace  des fonctions cuspidales, c'est-\`a-dire des fonctions $f\in C_{inv}(\mathfrak{g}({\mathbb F}_{q}))$ qui v\'erifient la condition suivante: soient $P$ un sous-groupe parabolique propre de $G$ et $M$ une composante de Levi de $P$; alors, pour tout $X\in \mathfrak{m}({\mathbb F}_{q})$, on a l'\'egalit\'e
$$\sum_{Y\in \mathfrak{u}_{P}({\mathbb F}_{q})}f(X+Y)=0.$$
Fixons une forme bilin\'eaire sym\'etrique non d\'eg\'en\'er\'ee $<.,.>$ sur $\mathfrak{g}({\mathbb F}_{q})$ invariante par conjugaison par $G({\mathbb F}_{q})$. Fixons aussi un caract\`ere non trivial $\psi:{\mathbb F}_{q}\to {\mathbb C}^{\times}$. On d\'efinit la transformation de Fourier $f\mapsto \hat{f}$ dans $C(\mathfrak{g}({\mathbb F}_{q}))$ par
$$\hat{f}(X)=q^{-dim(\mathfrak{g})/2}\sum_{Y\in \mathfrak{g}({\mathbb F}_{q})}f(Y)\psi(<X,Y>)$$
pour tout $X\in \mathfrak{g}({\mathbb F}_{q})$. 

On note $C_{nil}(\mathfrak{g}({\mathbb F}_{q}))$ le sous-espace de $C_{inv}(\mathfrak{g}({\mathbb F}_{q}))$ form\'e des fonctions \`a support nilpotent. On pose $C_{nil,cusp}(\mathfrak{g}({\mathbb F}_{q}))= C_{nil}(\mathfrak{g}({\mathbb F}_{q}))\cap C_{cusp}(\mathfrak{g}({\mathbb F}_{q}))$. On note $FC(\mathfrak{g}({\mathbb F}_{q}))$ le sous-espace des fonctions $f\in C_{nil}(\mathfrak{g}({\mathbb F}_{q}))$ telles que $\hat{f}$ appartient elle-aussi \`a $C_{nil}(\mathfrak{g}({\mathbb F}_{q}))$. Par d\'efinition,  l'espace $FC(\mathfrak{g}({\mathbb F}_{q}))$ est invariant par transformation de Fourier. 

On a 

(2) si $Z(G)^0\not=\{1\}$,  $FC(\mathfrak{g}({\mathbb F}_{q}))=0$. 

En effet, on a $FC(\mathfrak{g}({\mathbb F}_{q}))\subset FC^{Z(G)^0}(\mathfrak{z}(G)({\mathbb F}_{q}))\otimes_{{\mathbb C}} FC^{G_{SC}}(\mathfrak{g}_{SC}({\mathbb F}_{q}))$   et il suffit de prouver que $FC^{Z(G)^0}(\mathfrak{z}(G)({\mathbb F}_{q}))=\{0\}$.  L'espace des fonctions \`a support nilpotent dans $\mathfrak{z}(G)({\mathbb F}_{q})$ est la droite port\'ee par la fonction caract\'eristique de $\{0\}$. La transform\'ee de Fourier de celle-ci est une fonction constante non nulle. Cette derni\`ere  est \`a support nilpotent si et seulement si $\mathfrak{z}(G)({\mathbb F}_{q})=\{0\}$.  Cette condition \'equivaut \`a $Z(G)^0=\{1\}$. $\square$

Supposons donc $Z(G)^{0}=\{1\}$, c'est-\`a-dire $G$ semi-simple. Lusztig a prouv\'e en \cite{L} que $FC(\mathfrak{g}({\mathbb F}_{q}))$ avait pour base les fonctions caract\'eristiques des faisceaux-caract\`eres cuspidaux qui sont invariants par $\Gamma_{k}$, consid\'er\'ees \`a homoth\'etie pr\`es (l\`a encore, la notion de faisceau-caract\`ere d\'epend du groupe $G$ et pas seulement de $\mathfrak{g}$). En particulier, $FC(\mathfrak{g}({\mathbb F}_{q}))\subset C_{nil, cusp}(\mathfrak{g}({\mathbb F}_{q}))$.

\section{Groupes sur un corps $p$-adique}\label{groupespadiques}
Soit $p$ un nombre premier et $F$ une extension finie de ${\mathbb Q}_{p}$. On utilise les notations usuelles: $\mathfrak{o}_{F}$ est l'anneau des entiers de $F$, $\mathfrak{p}_{F}$ est son id\'eal maximal, ${\mathbb F}_{q}$ est le corps r\'esiduel, $\vert .\vert _{F}$ et $val_{F}$ sont les valeur absolue et valuation usuelles et on fixe une uniformisante $\varpi_{F}$. On note $F^{nr}$ la plus grande extension de $F$ contenue dans $\bar{F}$ et non ramifi\'ee sur $F$. 

Soit $G$ un groupe r\'eductif connexe d\'efini sur $F$. Notons  $rg(G)$ le rang de $G$ sur $\bar{F}$. On suppose

(1) $p\geq sup(2h(G)+1,rg(G)+2)$.

Cette hypoth\`ese implique que $G$ est d\'eploy\'e sur une extension galoisienne finie $F'$ de $F$ de degr\'e premier \`a $p$, a fortiori mod\'er\'ement ramifi\'ee. On note $A_{G}$ le plus grand tore  d\'eploy\'e contenu dans le centre de $G$,  $G_{reg}$ l'ensemble des \'el\'ements semi-simples fortement r\'eguliers de $G$ et $G_{ell}(F)$  le sous-ensemble de $G_{reg}(F)$ form\'e des \'el\'ements elliptiques, c'est-\`a-dire des $x\in G_{reg}(F)$ dont le centralisateur $T=Z_{G}(x)$ v\'erifie $A_{T}=A_{G}$ (la notation est contestable: il n'y a pas de sous-ensemble alg\'ebrique $G_{ell}$).  On utilise les notations analogues $\mathfrak{g}_{reg}$, $\mathfrak{g}_{ell}(F)$, pour l'alg\`ebre de Lie. 

On note $Imm(G_{AD})$ l'immeuble de Bruhat-Tits du groupe $G_{AD}$ sur $F$. Le groupe $G_{AD}(F)$, a fortiori le groupe $G(F)$, agit sur $Imm(G_{AD})$. Cet immeuble est muni d'une d\'ecomposition en facettes. Pour chaque facette ${\cal F}$, on introduit le groupe $K^{\dag}_{{\cal F}}$ des \'el\'ements de $G(F)$ dont l'action conserve la facette, le sous-groupe parahorique $K^0_{{\cal F}}$ et son radical pro-$p$-unipotent $K^+_{{\cal F}}$. Bruhat et Tits ont d\'efini un groupe r\'eductif connexe $G_{{\cal F}}$ d\'efini sur ${\mathbb F}_{q}$ tel que $K^0_{{\cal F}}/K^+_{{\cal F}}\simeq G_{{\cal F}}({\mathbb F}_{q})$. Aux groupes $K^0_{{\cal F}}$ et $K^+_{{\cal F}}$ sont associ\'ees des sous-$\mathfrak{o}_{F}$-alg\`ebres $\mathfrak{k}_{{\cal F},0}$ et $\mathfrak{k}_{{\cal F},0+}$ de $\mathfrak{g}(F)$. On a $\mathfrak{k}_{{\cal F},0}/\mathfrak{k}_{{\cal F},0+}\simeq \mathfrak{g}_{{\cal F}}({\mathbb F}_{q})$. Sous l'hypoth\`ese (1), on peut munir et on munit $\mathfrak{g}(F)$ de la mesure de Haar telle que $mes(\mathfrak{k}_{{\cal F},0})=\vert \mathfrak{g}_{{\cal F}}({\mathbb F}_{q})\vert ^{1/2}$ pour toute facette ${\cal F}$.  On a alors $mes(\mathfrak{k}_{{\cal F},0+})=\vert \mathfrak{g}_{{\cal F}}({\mathbb F}_{q})\vert ^{-1/2}$. On peut munir et on munit $G(F)$ de la mesure de Haar telle que $mes(K_{{\cal F}}^+)=mes(\mathfrak{k}_{{\cal F},0+})$ pour toute facette ${\cal F}$. Les m\^emes choix de mesures vaudront pour tout sous-groupe r\'eductif connexe de $G$.  Sous l'hypoth\`ese (1), on sait que l'on peut choisir et on choisit une forme bilin\'eaire sym\'etrique non d\'eg\'en\'er\'ee $<.,.>$ sur $\mathfrak{g}(F)$, invariante par l'action adjointe de $G_{AD}(F)$ et v\'erifiant la condition suivante. Pour tout $\mathfrak{o}_{F}$-r\'eseau $\mathfrak{h}\subset \mathfrak{g}(F)$, notons $\mathfrak{h}^{\perp}=\{X\in \mathfrak{g}(F); \forall Y\in \mathfrak{h}, <Y,X>\in \mathfrak{p}_{F}\}$. Alors $\mathfrak{k}_{{\cal F},0}^{\perp}=\mathfrak{k}_{{\cal F},0+}$ pour toute facette ${\cal F}$. On fixe un caract\`ere $\psi:F\to {\mathbb C}^{\times}$ de conducteur $\mathfrak{p}_{F}$. On d\'efinit la transformation de Fourier $f\mapsto \hat{f}$ dans $C_{c}^{\infty}(\mathfrak{g}(F))$ par
$$\hat{f}(X)=\int_{\mathfrak{g}(F)}f(Y)\psi(<Y,X>)\,dY.$$
La mesure que l'on a fix\'ee est autoduale pour cette transformation. Soit ${\cal F}$ une facette de l'immeuble. De la forme $<.,.>$ se d\'eduit naturellement une forme encore not\'ee $<.,.>$ sur $\mathfrak{g}_{{\cal F}}({\mathbb F}_{q})$. Du caract\`ere $\psi$ se d\'eduit aussi un caract\`ere encore not\'e $\psi$ de ${\mathbb F}_{q}$. Ces donn\'ees permettent de d\'efinir une transformation de Fourier dans $C(\mathfrak{g}_{{\cal F}}({\mathbb F}_{q}))$. Soit $f\in \mathfrak{g}_{{\cal F}}({\mathbb F}_{q})$. On rel\`eve $f$ en une fonction sur $\mathfrak{k}_{{\cal F},0}$ invariante par translations par $\mathfrak{k}_{{\cal F},0+}$. On prolonge cette fonction en une fonction sur $\mathfrak{g}(F)$ nulle hors de $\mathfrak{k}_{{\cal F},0}$. On note $f_{{\cal F}}$ la fonction ainsi obtenue. On v\'erifie que $(f_{{\cal F}})\hat{}=(\hat{f})_{{\cal F}}$. 

Pour $X\in \mathfrak{g}_{reg}(F)$ et $f\in C_{c}^{\infty}(\mathfrak{g}(F))$, on d\'efinit l'int\'egrale orbitale
$$I^G(X,f)=d^G(X)^{1/2}\int_{Z_{G}(X)(F)\backslash G(F)}f(g^{-1}Xg)\,dg.$$
Puisque $X$ est r\'egulier, le centralisateur $Z_{G}(X)$ est   un tore et il est muni  de la mesure analogue \`a celle que l'on a fix\'ee sur $G(F)$. Le terme $d^G(X)$ est l'habituel discriminant de Weyl. On \'etend cette d\'efinition \`a tout $X\in \mathfrak{g}(F)$ en posant simplement
$$I^G(X,f)=\int_{Z_{G}(X)(F)\backslash G(F)}f(g^{-1}Xg)\,dg$$
pour $X\in \mathfrak{g}(F)-\mathfrak{g}_{reg}(F)$, o\`u on fixe arbitrairement une mesure de Haar sur $Z_{G}(X)(F)$. 

   On note $I(\mathfrak{g}(F))$ le quotient de $C_{c}^{\infty}(\mathfrak{g}(F))$ par le sous-espace des fonctions $f$ telles que $I^G(X,f)=0$ pour tout $X\in \mathfrak{g}(F)$ (ou tout $X\in \mathfrak{g}_{reg}(F)$, les deux conditions sont \'equivalentes). Le dual $I(\mathfrak{g}(F))^*$ est l'espace des distributions sur $\mathfrak{g}(F)$ invariantes par conjugaison par $G(F)$. Les int\'egrales orbitales peuvent \^etre consid\'er\'ees comme des \'el\'ements de ce dual. 
On note $I_{cusp}(\mathfrak{g}(F))$ le sous-espace de $I(\mathfrak{g}(F))$ form\'e des \'el\'ements $f\in I(\mathfrak{g}(F))$ tels que $I^G(X,f)=0$ pour tout $X\in G_{reg}(F)-G_{ell}(F)$. La transformation de Fourier se descend en une transformation de $I(\mathfrak{g}(F))$ et celle-ci pr\'eserve le sous-espace $I_{cusp}(\mathfrak{g}(F))$. 
 
Pour $f,f'\in I_{cusp}(\mathfrak{g}(F))$, on d\'efinit le produit scalaire elliptique $J^G_{ell}(f',f)$ de la fa\c{c}on suivante. Fixons un ensemble de repr\'esentants  ${\cal T}_{ell}$ des classes de conjugaison par $G(F)$ des sous-tores maximaux elliptiques de $G$, c'est-\`a-dire des sous-tores maximaux $T$ tels que $A_{T}=A_{G}$. Pour tout tel tore $T$, posons $W^G(T)=Norm_{G(F)}(T)/T(F)$. Alors
$$(2) \qquad J^G_{ell}(f',f)=\sum_{T\in {\cal T}_{ell}}\vert W^G(T)\vert ^{-1}mes(A_{G}(F)\backslash T(F))\int_{\mathfrak{t}(F)}I^G(X,\bar{f}')I^G(X,f)\,dX.$$
Ce produit est d\'efini positif. 

{\bf Remarque.} Conform\'ement \`a nos conventions, la mesure sur $A_{G}(F)$ est  ici la mesure analogue \`a celle que l'on a fix\'ee sur $G(F)$. Dans ses articles, Arthur  choisit celle pour laquelle  la mesure du plus sous-groupe compact maximal de $A_{G}(F)$ vaut $1$. Cela ne change rien aux propri\'et\'es que nous utiliserons du produit scalaire.
\bigskip

Les donn\'ees $F$ et $G$ seront d\'esormais conserv\'ees pour tout l'article.

\section{R\'eseaux de Moy-Prasad}\label{reseaux}

Soit $x\in Imm(G_{AD})$. Pour tout  r\'eel $r\geq0$, Moy et Prasad ont d\'efini un sous-groupe $K_{x,r}\subset G(F)$. Si $s>r$, $G_{x,s}$ est un sous-groupe distingu\'e de $G_{x,r}$ et l'application $r\mapsto G_{x,r}$ est localement constante \`a gauche. On pose $G_{x,r+}=\cup_{s>r}G_{x,s}=G_{x,r+\epsilon}$ pour $\epsilon>0$ assez petit.  En notant ${\cal F}$ la facette \`a laquelle appartient $x$, on a $K_{x,0}=K_{{\cal F}}^0$ et $K_{x,0+}=K_{{\cal F}}^+$. Aux groupes $G_{x,r}$ et $G_{x,r+}$ sont associ\'ees des sous-$\mathfrak{o}_{F}$-alg\`ebres $\mathfrak{g}_{x,r}$ ou $\mathfrak{g}_{x,r+}$ de $\mathfrak{g}(F)$. On peut d'ailleurs d\'efinir ces objets pour tout $r\in {\mathbb R}$, sans la restriction $r\geq0$ (pour $r<0$, ce ne sont plus des alg\`ebres, seulement des $\mathfrak{o}_{F}$-r\'eseaux). On  a les \'egalit\'es $\mathfrak{g}_{x,r+1}=\mathfrak{p}_{F}\mathfrak{g}_{x,r}$ et $\mathfrak{g}_{x,r}^{\perp}=\mathfrak{g}_{x,(-r)+}$. 

On note $\mathfrak{g}(F)_{r}$ la r\'eunion des $\mathfrak{g}_{x,r}$ quand $x$ parcourt $Imm(G_{AD})$. On d\'efinit $\mathfrak{g}(F)_{r+}$ de fa\c{c}on similaire. Dans le cas $r=0$, on pose  $\mathfrak{g}_{ent}(F)=\mathfrak{g}(F)_{0}$, $\mathfrak{g}_{tn}(F)=\mathfrak{g}(F)_{0+}$. Les \'el\'ements de $\mathfrak{g}_{ent}(F)$ sont dits entiers, ceux de $\mathfrak{g}_{tn}(F)$ sont dit topologiquement nilpotents.

Pr\'ecisons nos notations en y glissant le corps de base $F$: $Imm_{F}(G_{AD})$, $\mathfrak{g}_{F,x,r}$. Fixons une extension galoisienne finie $F'/F$ de degr\'e premier \`a $p$, telle que $G$ soit d\'eploy\'e sur $F'$. Notons $e(F'/F)$ l'indice de ramification. Le groupe $\Gamma_{F'/F}$ agit sur $Imm_{F'}(G_{AD})$ et $Imm_{F}(G_{AD})$ s'identifie canoniquement au sous-ensemble des points fixes $(Imm_{F'}(G_{AD}))^{\Gamma_{F'/F}}$. L'action de $G(F)$ est la restriction de celle de $G(F')$. Pour tous $x\in Imm_{F}(G_{AD})$ et $r\in {\mathbb R}$, on a $\mathfrak{g}_{F,x,r}=\mathfrak{g}_{F',x,er}\cap \mathfrak{g}(F)$. Cela ram\`ene la description de ces  r\'eseaux au cas o\`u $G$ est d\'eploy\'e.

Supposons donc $G$ d\'eploy\'e sur $F$. Fixons un sous-tore maximal d\'eploy\'e $T$ de $G$.  Pour tout $r\in {\mathbb R}$, on note $\mathfrak{t}(F)_{r}$ le sous-ensemble des $X\in \mathfrak{t}(F)$ tels que $val_{F}(x^*(X))\geq r$ pour tout $x^*\in X^*(T)$. Au tore $T$ est associ\'e un appartement dans $Imm(G_{AD})$ qui est un espace affine sous ${\cal T}/{\cal A}_{G}$.  Fixons un sommet hypersp\'ecial $s$ de cet appartement. On identifie celui-ci \`a ${\cal T}/{\cal A}_{G}$ en envoyant $s$ sur $0$. Notons $\Sigma$ l'ensemble des racines de $T$ dans $G$ et, pour tout $\alpha\in \Sigma$, notons $\mathfrak{u}_{\alpha}\subset \mathfrak{g}$ la droite radicielle associ\'ee \`a $\alpha$.  On peut en fixer un vecteur de base $E_{\alpha}$  de sorte que
$$\mathfrak{g}_{s,0}=\mathfrak{t}_{0}\oplus(\oplus_{\alpha\in \Sigma} \mathfrak{o}_{F}E_{\alpha}).$$
Soient  $x\in    {\cal T}/{\cal A}_{G}$ et $r\in {\mathbb R}$. Pour $\alpha\in \Sigma$, notons $n_{x,r,\alpha}$, resp. $n_{x,r+,\alpha}$,  le plus petit entier relatif $n$ tel que $r-\alpha(x)\leq n$, resp. $r-\alpha(x)<n$. Alors
$$\mathfrak{g}_{x,r}=\mathfrak{t}(F)_{r}\oplus(\oplus_{\alpha\in \Sigma} \mathfrak{p}^{n_{x,r,\alpha}}_{F}E_{\alpha}),$$
$$\mathfrak{g}_{x,r+}=\mathfrak{t}(F)_{r+}\oplus(\oplus_{\alpha\in \Sigma} \mathfrak{p}^{n_{x,r+,\alpha}}_{F}E_{\alpha}).$$
 
\section{Quelques espaces de fonctions}\label{quelquesespaces}

Pour un sous-$\mathfrak{o}_{F}$-module $\mathfrak{h}\subset \mathfrak{g}(F)$ et un sous-$\mathfrak{o}_{F}$-r\'eseau $\mathfrak{l}\subset \mathfrak{h}$, notons $C_{c}(\mathfrak{h}/\mathfrak{l})$ 
 le sous-espace de $C_{c}^{\infty}(\mathfrak{g}(F))$ form\'e des fonctions \`a support dans $\mathfrak{h}$ et  invariantes par translations par $\mathfrak{l}$. On supprime l'indice $c$ si $\mathfrak{h}$ est compact.  Posons
$${\cal H}_{r}=\sum_{x\in Imm(G_{AD})}C_{c}(\mathfrak{g}(F)/\mathfrak{g}_{x,r}),$$
$${\cal H}_{r+}=\sum_{x\in Imm(G_{AD})}C_{c}(\mathfrak{g}(F)/\mathfrak{g}_{x,r+}),$$
$${\cal H}_{r,r+}=\sum_{x\in Imm(G_{AD})}C(\mathfrak{g}_{x,r}/\mathfrak{g}_{x,r+}).$$
On a \'evidemment ${\cal H}_{r,r+}\subset {\cal H}_{r+}\subset {\cal H}_{s}$ pour tout $s>r$.

Soit $x\in Imm(G_{AD})$. Notons $\mathfrak{ L}_{x,r}$ l'ensemble des $\mathfrak{o}_{F}$-r\'eseaux $\mathfrak{l}\subset \mathfrak{g}(F)$ tels que, pour tout voisinage $V$ de $x$ dans $Imm(G_{AD})$, il existe $y\in V$ de sorte que $\mathfrak{l}=\mathfrak{g}_{y,r}$.   On montrera ci-dessous que

(1) on a $\mathfrak{g}_{x,r+}\subset \mathfrak{l}\subset \mathfrak{g}_{x,r}$ pour tout $\mathfrak{l}\in \mathfrak{ L}_{x,r}$.
 
 A fortiori, $\mathfrak{L}_{x,r}$ est fini. On pose
 $${\cal H}_{r}^{\sharp}=\sum_{x\in Imm(G_{AD})}\sum_{\mathfrak{l}\in \mathfrak{L}_{x,r}}C(\mathfrak{g}_{x,r}/\mathfrak{l}).$$
  Pour $x$ et $\mathfrak{l}$ intervenant dans cette somme, il existe par d\'efinition $y\in Imm(G_{AD})$ tel que $\mathfrak{l}=\mathfrak{g}_{y,r}$ donc $C(\mathfrak{g}_{x,r}/\mathfrak{l})\subset C_{c}(\mathfrak{g}(F)/\mathfrak{g}_{y,r})$. D'o\`u l'inclusion ${\cal H}_{r}^{\sharp}\subset {\cal H}_{r}$. 
  
  \begin{lem}{Soit $r\in {\mathbb R}$. Il existe un r\'eel $\epsilon_{0}>0$ tel que
  
  (i) on a les \'egalit\'es $\mathfrak{g}(F)_{r-\epsilon}=\mathfrak{g }(F)_{r}$ et ${\cal H}_{r-\epsilon}={\cal H}_{r}$ pour tout $\epsilon\in [0,\epsilon_{0}[$;
  
  (ii) l'espace ${\cal H}_{r-\epsilon,(r-\epsilon)+}$ ne d\'epend pas de $\epsilon$ pour $\epsilon\in ]0,\epsilon_{0}[$ et est inclus dans ${\cal H}_{r}^{\sharp}$.}\end{lem}
  
  {\bf Remarque.} La premi\`ere propri\'et\'e de (i) r\'esulte  du fait connu que l'application $r\mapsto \mathfrak{g}(F)_{r}$ est localement constante \`a gauche et que ses sauts sont un sous-ensemble discret de ${\mathbb R}$. La deuxi\`eme propri\'et\'e de (i) nous semble moins \'evidente. 
  \bigskip
  
  Preuve. On fixe $r\in {\mathbb R}$. L'immeuble $Imm(G_{AD})$ est muni d'une distance, on note $\vert y-x\vert $ la distance entre deux points $x,y$ de l'immeuble. Montrons que
  
  (2) pour tout $x\in Imm(G_{AD})$, il existe $\epsilon_{x},\eta_{x}>0$ tels que, pour tout $\epsilon\in [0,\epsilon_{x}[$ et pour tout $y\in Imm(G_{AD})$ tel que $\vert y-x\vert <\eta_{x}$, on a $\mathfrak{g}_{x,r+}\subset \mathfrak{g}_{y,(r-\epsilon)+}\subset \mathfrak{g}_{y,r-\epsilon}\subset \mathfrak{g}_{x,r}$.
  
  Fixons une extension galoisienne finie $F'/F$ de degr\'e premier \`a $p$ telle que $G$ soit d\'eploy\'e sur $F'$. Si l'on d\'emontre l'assertion (2) pour le corps de base $F'$, la m\^eme assertion s'en d\'eduit sur le corps de base $F$ en prenant les invariants par $\Gamma_{F'/F}$ de chacun des r\'eseaux. Autrement dit, en oubliant cette extension $F'$, on peut supposer $G$ d\'eploy\'e sur $F$. Un voisinage de $x$ dans l'immeuble est contenu dans  la r\'eunion d'un nombre fini d'appartements. On peut fixer un appartement    associ\'e \`a un tore maximal d\'eploy\'e $T$ et auquel appartient $x$  et se limiter \`a consid\'erer des points $y$ dans cet appartement. Avec les notations du paragraphe pr\'ec\'edent, l'assertion (2) se traduit par les deux assertions suivantes:
  
  (3) il existe $\epsilon_{1}>0$ tel que $\mathfrak{t}(F)_{r-\epsilon}=\mathfrak{t}(F)_{r}$ pour tout $\epsilon\in [0,\epsilon_{1}[$;
  
  (4) il existe $\epsilon_{x},\eta_{x}>0$ tels que, pour tout $\epsilon\in [0,\epsilon_{x}[$, pour tout $y\in {\cal T}/{\cal A}_{G}$ tel que $\vert y-x\vert <\eta_{x}$ et pour tout $\alpha\in \Sigma$, on a les in\'egalit\'es $n_{x,r,\alpha}\leq n_{y,r-\epsilon,\alpha}\leq n_{y,(r-\epsilon)+,\alpha}\leq n_{x,r+,\alpha}$.
  
 L'assertion (3) r\'esulte imm\'ediatement de la d\'efinition de $\mathfrak{t}(F)_{r}$. Pour tout $\alpha\in \Sigma$, on a les in\'egalit\'es $n_{x,r,\alpha}-1<r-\alpha(x)\leq n_{x,r,\alpha}$ et $n_{x,r+,\alpha}-1\leq r-\alpha(x)<n_{x,r+,\alpha}$.  On peut fixer $\eta_{x}>0$ de sorte que, pour $y\in {\cal T}/{\cal A}_{G}$, la condition $\vert y-x\vert <\eta_{x}$ entra\^{\i}ne $\alpha(y)-\alpha(x)<\frac{1}{2}(r-\alpha(x)-n_{x,r,\alpha}+1)$  et $\alpha(x)-\alpha(y)<n_{x,r+,\alpha}-r+\alpha(x)$  pour tout $\alpha\in \Sigma$. On peut fixer $\epsilon_{x}>0$ tel que $\epsilon_{x}< \frac{1}{2}(r-\alpha(x)-n_{x,r,\alpha}+1)$  pour tout $\alpha\in \Sigma$. Supposons alors $\epsilon\in [0,\epsilon_{x}[$ et $\vert y-x\vert <\eta_{x}$. Soit $\alpha\in \Sigma$. On a 
 $$r-\epsilon-\alpha(y)= r-\alpha(x)-\epsilon+\alpha(x)-\alpha(y)> r-\alpha(x)-(r-\alpha(x)-n_{x,r,\alpha}+1)=n_{x,r,\alpha}-1.$$
 Par d\'efinition, cela entra\^{\i}ne $n_{y,r-\epsilon,\alpha}\geq n_{x,r,\alpha}$. On a aussi
 $$r-\epsilon-\alpha(y)\leq r-\alpha(y)=r-\alpha(x)-\alpha(y)+\alpha(x)<r-\alpha(x)+(n_{x,r+,\alpha}-r+\alpha(x))=n_{x,r+,\alpha}.$$
 Par d\'efinition, cela entra\^{\i}ne $n_{y,(r-\epsilon)+,\alpha}\leq n_{x,r+,\alpha}$. Cela d\'emontre (4), d'o\`u (2). 
 
 L'assertion (2) entra\^{\i}ne l'assertion (1) pr\'ec\'edant l'\'enonc\'e. En fait, pour $x\in Imm(G_{AD})$ et $\eta_{x}$ comme en (2), on a
 
 (5) $\mathfrak{L}_{x,r}$ est l'ensemble des r\'eseaux $\mathfrak{g}_{y,r}$ pour $y\in Imm(G_{AD})$ tel que $\vert y-x\vert <\eta_{x}$.
 
 En effet,   $\mathfrak{L}_{x,r}$ est contenu dans cet ensemble de r\'eseaux par d\'efinition. Inversement, soit $y\in Imm(G_{AD})$ tel que $\vert y-x\vert <\eta_{x}$. Notons $[x,y]$ la g\'eod\'esique joignant $x$ \`a $y$, que l'on peut identifier au segment r\'eel  $[0,1]$. Pour $t\in ]0,1]$, notons $y_{t}$ l'\'el\'ement de la g\'eod\'esique qui s'identifie \`a $t\in [0,1]$. 
 Comme dans la preuve de (2), on \'etend le corps de d\'efinition $F$ et on fixe un appartement associ\'e \`a un tore maximal d\'eploy\'e $T$ qui contient $x$ et $y$.  Un calcul analogue \`a ceux de la preuve de (2) montre que $\mathfrak{g}_{y_{t},r}=\mathfrak{g}_{y,r}$. Pour tout voisinage $V$ de $x$ dans $Imm(G_{AD})$, il existe $t\in ]0,1]$ tel que $y_{t}\in V$. En appliquant la d\'efinition de $\mathfrak{L}_{x,r}$, il en r\'esulte que $\mathfrak{g}_{y,r}\in \mathfrak{L}_{x,r}$. D'o\`u (5).
 
 Il existe un sous-ensemble compact $\Delta\subset Imm(G_{AD})$ tel que  $Imm(G_{AD})=\cup_{g\in G(F)}g\Delta$. Fixons un tel $\Delta$. Pour tout $x\in \Delta$, fixons $\epsilon_{x}$ et $\eta_{x}$ de sorte que (2) soit v\'erifi\'e. L'ensemble $\Delta$ est recouvert par  les ouverts $V_{x}=\{y\in Imm(G_{AD}); \vert y-x\vert <\eta_{x}\}$ quand $x$ parcourt $\Delta$. On en extrait un recouvrement fini $\Delta\subset \cup_{i=1,...,n}V_{x_{i}}$. Posons $\epsilon_{2}=inf_{i=1,...,n}\epsilon_{x_{i}}$. Notons ${\cal R}_{r,\epsilon_{2}}$ l'ensemble des r\'eseaux $\mathfrak{g}_{x,r-\epsilon}$ quand $x$ d\'ecrit $Imm(G_{AD})$ et $\epsilon$ d\'ecrit $[0,\epsilon_{2}[$. Montrons que
 
 (6) l'ensemble des orbites de l'action de $G(F)$ dans ${\cal R}_{r,\epsilon_{2}}$ est fini.
 
 On peut se limiter \`a consid\'erer les    r\'eseaux $\mathfrak{g}_{x,r-\epsilon}$ quand $x$ d\'ecrit $\Delta$ et $\epsilon$ d\'ecrit $[0,\epsilon_{2}[$. Or, d'apr\`es (2),  ces r\'eseaux appartiennent \`a la r\'eunion sur $i=1,...,n$ des ensembles de $\mathfrak{o}_{F}$ r\'eseaux $\mathfrak{l}$ tels que $\mathfrak{g}_{x_{i},r+}\subset \mathfrak{l}\subset \mathfrak{g}_{x_{i}}$. Ces derniers ensembles \'etant finis, (6) s'ensuit. 
 
 Notons ${\cal R}_{r}$ l'ensemble des r\'eseaux $\mathfrak{g}_{x,r}$ quand $x$ d\'ecrit $Imm(G_{AD})$. Montrons que
 
 (7) il existe $\epsilon_{3}>0$ tel que, pour $\epsilon\in [0,\epsilon_{3}[$, on a ${\cal R}_{r}={\cal R}_{r-\epsilon}$.
 
 D'apr\`es (6), on peut fixer des points $z_{j}\in Imm(G_{AD})$, $j=1,...,m$ de sorte que, pour tout $x\in Imm(G_{AD})$ et tout $\epsilon\in [0,\epsilon_{2}[$, il existe $j\in \{1,...,m\}$ et $g\in G(F)$ de sorte que $\mathfrak{g}_{x,r-\epsilon}=Ad(g)(\mathfrak{g}_{z_{j},r-\epsilon})$. Pour tout $j$, il existe $\epsilon_{3,j}>0$ de sorte que $\mathfrak{g}_{z_{j},r-\epsilon}=\mathfrak{g}_{z_{j},r}$ pour tout $\epsilon\in [0,\epsilon_{3,j}[$. On pose $\epsilon_{3}=inf(\epsilon_{2},inf_{j}\epsilon_{3,j})$. 
   Pour $x\in Imm(G_{AD})$ et $\epsilon,\epsilon'\in [0,\epsilon_{3}[$, soit $g\in G(F)$ et $j\in \{1,...,m\}$ de sorte que  $\mathfrak{g}_{x,r-\epsilon}=Ad(g)(\mathfrak{g}_{z_{j},r-\epsilon})$. On a aussi $\mathfrak{g}_{x,r-\epsilon}=Ad(g)(\mathfrak{g}_{z_{j},r-\epsilon'})=\mathfrak{g}_{gz_{j},r-\epsilon'}$, donc $\mathfrak{g}_{x,r-\epsilon}$ appartient \`a ${\cal R}_{r-\epsilon'}$. Cela prouve (7).
   
   L'assertion (7) et les d\'efinitions entra\^{\i}nent que $\mathfrak{g}(F)_{r}= \mathfrak{g}(F)_{r-\epsilon}$ et ${\cal H}_{r}={\cal H}_{r-\epsilon}$ pour $\epsilon\in [0,\epsilon_{3}[$. Cela d\'emontre  le (i) de l'\'enonc\'e (en supposant $\epsilon_{0}\leq \epsilon_{3}$). 
   
   Montrons que
   
   (8) pour $\epsilon,\epsilon'$ tels que $0<\epsilon'<\epsilon<\epsilon_{2}$, on a ${\cal H}_{r-\epsilon,(r-\epsilon)+}\subset {\cal H}_{r-\epsilon',(r-\epsilon')+}$.
   
   Avec les notations  d\'ej\`a introduites, on peut fixer $i\in \{1,...,n\}$, $y\in V_{x_{i}}$ et il suffit de prouver qu'il existe $z\in  [x_{i},y]$ tel que $\mathfrak{g}_{z,r-\epsilon'}=\mathfrak{g}_{y,r-\epsilon}$ et $\mathfrak{g}_{z,(r-\epsilon')+}=\mathfrak{g}_{y,(r-\epsilon)+}$.  On pose simplement $x=x_{i}$. Comme dans la preuve de (2), on \'etend le corps de d\'efinition $F$ et on fixe un appartement associ\'e \`a un tore maximal d\'eploy\'e $T$ qui contient $x$ et $y$. On voit que la condition $\mathfrak{g}_{z,r-\epsilon'}=\mathfrak{g}_{y,r-\epsilon}$ et $\mathfrak{g}_{z,(r-\epsilon')+}=\mathfrak{g}_{y,(r-\epsilon)+}$ \'equivaut \`a la r\'eunion sur $\alpha\in \Sigma$ des conditions suivantes
   
   $ r-\epsilon-\alpha(y)<n_{x,r,\alpha}$, resp. $ r-\epsilon-\alpha(y)=n_{x,r,\alpha}$, resp. $ r-\epsilon-\alpha(y)>n_{x,r,\alpha}$, \'equivaut \`a $ r-\epsilon'-\alpha(z)<n_{x,r,\alpha}$, resp. $r-\epsilon'-\alpha(z)=n_{x,r,\alpha}$, resp. $r-\epsilon'-\alpha(z)>n_{x,r,\alpha}$.
   
   Pour une racine $\alpha$ telle que $r-\alpha(x)<n_{x,r,\alpha}$, on a  $ r-\epsilon-\alpha(y)<n_{x,r,\alpha}$ et $ r-\epsilon'-\alpha(z)<n_{x,r,\alpha}$ d'apr\`es la d\'efinition de $\eta_{x_{i}}$ et $\epsilon_{x_{i}}$. La condition ci-dessus est donc automatique. Supposons $r-\alpha(x)=n_{x,r,\alpha}$. Alors les conditions  $ r-\epsilon-\alpha(y)<n_{x,r,\alpha}$, resp. $ r-\epsilon-\alpha(y)=n_{x,r,\alpha}$, resp. $ r-\epsilon-\alpha(y)>n_{x,r,\alpha}$, \'equivalent \`a $\alpha(x)-\alpha(y)<\epsilon$, resp. $\alpha(x)-\alpha(y)=\epsilon$, resp. $\alpha(x)-\alpha(y)>\epsilon$. On traduit de fa\c{c}on analogue les conditions concernant $z$ et $\epsilon'$. Pour $z=y_{\epsilon'/\epsilon}$,  on a $\alpha(x)-\alpha(z)=\frac{\epsilon'}{\epsilon}(\alpha(x)-\alpha(y))$ et on voit que  les conditions relatives \`a $y$ et $\epsilon$ sont bien \'equivalentes \`a celles relatives \`a $z$ et $\epsilon'$. Cela prouve (8).
   
   Il r\'esulte de (6), ou plus exactement de sa preuve, que, quand $\epsilon$ d\'ecrit $]0,\epsilon_{2}[$, les espaces ${\cal H}_{r-\epsilon,(r-\epsilon)+}$ ne parcourent qu'un nombre fini de sous-espaces de $C_{c}^{\infty}(\mathfrak{g}(F))$. L'assertion (8) dit que l'application $\epsilon\mapsto {\cal H}_{r-\epsilon,(r-\epsilon)+}$ est d\'ecroissante. Elle est donc constante au voisinage de $0$. C'est la premi\`ere assertion du (ii) de l'\'enonc\'e. Pour la seconde assertion, il suffit de prouver que, pour $\epsilon\in ]0,\epsilon_{2}[$ et pour $y\in Imm(G_{AD})$, il existe $x\in Imm(G_{AD})$ et $\mathfrak{l}\in \mathfrak{L}_{x,r}$ de sorte que $\mathfrak{l}\subset \mathfrak{g}_{y,(r-\epsilon)+}\subset \mathfrak{g}_{y,r-\epsilon}\subset \mathfrak{g}_{x,r}$. A conjugaison pr\`es, on peut supposer qu'il existe $i\in \{1,...,n\}$ tel que $y\in V_{x_{i}}$. Alors le point $x=x_{i}$ et le r\'eseau $\mathfrak{l}=\mathfrak{g}_{y,r}$ v\'erifient ces conditions d'apr\`es (2) et (5). $\square$
   
   On note ${\cal H}_{r}^{\natural}$ l'espace ${\cal H}_{r-\epsilon,(r-\epsilon)+}$ pour un $\epsilon\in ]0,\epsilon_{0}[$ quelconque. Il est inclus dans ${\cal H}_{r}^{\sharp}$.
   
Consid\'erons le cas $r=0$.   Soit $x\in Imm(G_{AD})$, notons ${\cal F}$ la facette \`a laquelle il appartient. On a $\mathfrak{g}_{x,0}=\mathfrak{k}_{{\cal F},0}$. Il r\'esulte de  (5) que, quand $\mathfrak{l}$ d\'ecrit $\mathfrak{L}_{x,0}$, les espaces $C(\mathfrak{g}_{x,0}/\mathfrak{l})$ s'identifient aux espaces de fonctions sur $\mathfrak{g}_{{\cal F}}({\mathbb F}_{q})/\mathfrak{p}({\mathbb F}_{q})$ o\`u $P$ d\'ecrit les sous-groupes paraboliques de $G_{{\cal F}}$. Ce sont aussi les images par transformation de Fourier des espaces de fonctions sur $\mathfrak{u}_{P}({\mathbb F}_{q})$. La somme sur $P$ de ces derniers espaces n'est autre que l'espace des fonctions sur $\mathfrak{g}_{{\cal F}}({\mathbb F}_{q})$ \`a support nilpotent.

\section{Le r\'esultat de DeBacker}\label{debacker}
 Notons $I(\mathfrak{g}(F))^*_{nil}
$ le sous-espace des distributions invariantes sur $\mathfrak{g}(F)$ \`a support nilpotent. Si on fixe un ensemble de repr\'esentants ${\cal N}$ des classes de conjugaison par $G(F)$ dans $\mathfrak{g}_{nil}(F)$, $I(\mathfrak{g}(F))^*_{nil}$
 a pour base les int\'egrales orbitales $I^G(N,.)$ pour $N\in {\cal N}$. Soit $r\in {\mathbb R}$. On note $I(\mathfrak{g}(F))^*_{r}
$, resp. $I(\mathfrak{g}(F))^*_{r+}$, le sous-espace des distributions invariantes sur $\mathfrak{g}(F)$ \`a support dans $\mathfrak{g}(F)_{r}$, resp. $\mathfrak{g}(F)_{r+}$. 
  On a les homomorphismes \'evidents
 
 $$I(\mathfrak{g}(F))^*_{nil}\stackrel{a_{r}}{\to}I(\mathfrak{g}(F))^*_{r}\stackrel{b_{r}}{\to} {\cal H}_{r}^*\stackrel{c_{r}}{\to}{\cal H}_{r}^{\sharp,*}\stackrel{d_{r}}{\to} {\cal H}_{r}^{\natural,*}.$$
 
 \begin{prop}{Les applications $b_{r}$ et $b_{r}\circ a_{r}$ ont m\^eme image. Les restrictions de $c_{r}$ et $d_{r}\circ c_{r}$ \`a l'image de $b_{r}$ sont injectives. Les applications $b_{r}\circ a_{r}$, $c_{r}\circ b_{r}\circ a_{r}$ et $d_{r}\circ c_{r}\circ b_{r}\circ a_{r}$ sont injectives.}\end{prop}
 
 Preuve.  Soit $\epsilon_{0}$ comme dans le lemme pr\'ec\'edent et $\epsilon\in ]0,\epsilon_{0}[$. On peut remplacer ci-dessus les espaces $I(\mathfrak{g}(F))^*_{r}$, $ {\cal H}_{r}^*$ et $ {\cal H}_{r}^{\natural,*}$ par $I(\mathfrak{g}(F))^*_{(r-\epsilon)+}$, ${\cal H}_{(r-\epsilon)+}^*$ et ${\cal H}^*_{r-\epsilon,(r-\epsilon)+}$. Alors le th\'eor\`eme 2.1.5 de \cite{D} dit que les
  applications $b_{r}$ et $b_{r}\circ a_{r}$ ont m\^eme image et que la restriction de $d_{r}\circ c_{r}$ \`a l'image de $b_{r}$ est injective. Evidemment, la restriction de $c_{r}$ \`a cette image l'est aussi. Pour $\lambda\in F^{\times}$ et $f\in C_{c}^{\infty}(\mathfrak{g}(F))$, notons $f^{\lambda}$ la fonction $f^{\lambda}(X)=f(\lambda X)$. On peut fixer une famille $(f_{N})_{N\in {\cal N}}$ d'\'el\'ements de $C_{c}^{\infty}(\mathfrak{g}(F))$ de sorte que  la matrice $(I^G(N,f_{N'}))_{N,N'\in {\cal N}}$ soit inversible. On sait que, pour tout $N\in {\cal N}$,, il existe un entier $d(N)\in {\mathbb Z}$ de sorte que $I^G(N,f^{\lambda^2})=\vert \lambda\vert _{F}^{d(N)}I^G(N,f)$ pour tout $\lambda\in F^{\times}$. La matrice $(I^G(N,f^{\lambda^2}_{N'}))_{N,N'\in {\cal N}}$ est donc elle-aussi inversible. Mais on peut choisir $\lambda$ tel que toutes les fonctions $f_{N}^{\lambda^2}$ appartiennent \`a ${\cal H}_{r}$. Il en r\'esulte que $b_{r}\circ a_{r}$ est injective. Alors $c_{r}\circ b_{r}\circ a_{r}$ et $d_{r}\circ c_{r}\circ b_{r}\circ a_{r}$ sont aussi injectives puisque les restrictions de $c_{r}$ et $d_{r}\circ c_{r}$ \`a l'image de $b_{r}$ sont injectives. $\square$
  
  \section{L'espace ${\cal D}(\mathfrak{g}(F))$ et une formule de produit}\label{produitscalaire}
Notons $S(G)$ l'ensemble des sommets de l'immeuble $Imm(G_{AD})$. Posons
$${\cal S}_{cusp}(\mathfrak{g}(F))=\sum_{s\in S(G)}C_{cusp}(\mathfrak{g}_{s}({\mathbb F}_{q})).$$
Ici, chaque espace $C_{cusp}(\mathfrak{g}_{s}({\mathbb F}_{q}))$ est identifi\'e \`a un sous-espace de $C_{c}^{\infty}(\mathfrak{g}(F))$ par l'application $f\mapsto f_{s}$. L'espace ${\cal S}_{cusp}(\mathfrak{g}(F))$ est donc lui-aussi  un sous-espace de $C_{c}^{\infty}(\mathfrak{g}(F))$. En fait, les \'el\'ements de ${\cal S}_{cusp}(\mathfrak{g}(F))$ sont cuspidaux, c'est-\`a-dire que l'image de ${\cal S}_{cusp}(\mathfrak{g}(F))$ dans $I(\mathfrak{g}(F))$ est contenue dans $I_{cusp}(\mathfrak{g}(F))$, cf. \cite{W2} lemme 10.1 (ce lemme concerne des fonctions sur le groupe $G(F)$ mais la m\^eme d\'emonstration s'applique \`a des fonctions sur $\mathfrak{g}(F)$).

 Pour $f\in {\cal S}_{cusp}(\mathfrak{g}(F))$, on d\'efinit une distribution invariante $D^G_{f}$ sur $\mathfrak{g}(F)$ par
$$D^G_{f}(f')=\int_{A_{G}(F)\backslash G(F)}\int_{\mathfrak{g}(F)}f'(g^{-1}Xg)\bar{f}(X)\,dX\,dg.$$
La double int\'egrale est convergente dans cet ordre. 

Soit $M$ un groupe de Levi de $G$. On d\'efinit l'espace ${\cal S}_{cusp}(\mathfrak{m}(F))$. Pour $f\in {\cal S}_{cusp}(\mathfrak{m}(F))$, on d\'efinit la distribution $D^M_{f}\in I(\mathfrak{m}(F))^*$ puis la distribution induite $D^G_{f}=Ind_{M}^G(D^M_{f})\in I(\mathfrak{g}(F))^*$. Posons
$${\cal S}(\mathfrak{g}(F))=\sum_{M}{\cal S}(\mathfrak{m}(F)),$$
o\`u $M$ parcourt les groupes de Levi de $G$. On a alors un homomorphisme antilin\'eaire $D^G:{\cal S}(\mathfrak{g}(F))\to I(\mathfrak{g}(F))^*$.  Le groupe $G(F)$ agit naturellement sur l'espace ${\cal S}(\mathfrak{g}(F))$: un \'el\'ement $g\in G(F)$ transporte un Levi $M$ sur le Levi $M'=gMg^{-1}$, un sommet $s\in S(M)$ sur un sommet $gs\in S(M')$ et une fonction $f\in C_{cusp}(\mathfrak{m}_{s}({\mathbb F}_{q}))$ sur une fonction $^gf\in C_{cusp}(\mathfrak{m}'_{gs}({\mathbb F}_{q}))$. On note ${\cal D}(\mathfrak{g}(F))$ l'espace des coinvariants de ${\cal S}(\mathfrak{g}(F))$ pour cette action de $G(F)$. L'homomorphisme $D^G$ se quotiente en un homomorphisme antilin\'eaire  encore not\'e $D^G:{\cal D}(\mathfrak{g}(F))\to I(\mathfrak{g}(F))^*$.

 Il est facile de d\'ecrire concr\`etement l'espace ${\cal D}(\mathfrak{g}(F))$. Soit $M$ un Levi de $G$. Le groupe $Norm_{G(F)}(M)$ agit naturellement sur $S(M)$. On fixe un ensemble de repr\'esentants $\underline{S}^G(M)$ pour cette action. Pour tout $s\in S(M)$, notons $K_{s}^{\dag,G}$ le sous-groupe des $g\in Norm_{G(F)}(M)$ tels que $gs=s$. Ce groupe agit naturellement sur $C_{cusp}(\mathfrak{m}_{s}({\mathbb F}_{q}))$ et on note comme toujours $C_{cusp}(\mathfrak{m}_{s}({\mathbb F}_{q}))^{K_{s}^{\dag,G}}$ le sous-espace des invariants. 
 Fixons un ensemble de repr\'esentants $\underline{{\cal L}}$ des classes de conjugaison par $G(F)$ de Levi de $G$.  On voit alors que le sous-espace
 $$\oplus_{M\in \underline{{\cal L}}}\oplus_{s\in \underline{S}^G(M)}  C_{cusp}(\mathfrak{m}_{s}({\mathbb F}_{q}))^{K_{s}^{\dag,G}}$$
 de ${\cal S}(\mathfrak{g}(F))$ s'envoie bijectivement sur ${\cal D}(\mathfrak{g}(F))$.
 
 Soient $M,M'$ deux Levi de $G$, $s\in S(M)$, $s'\in S(M')$, $f\in C_{cusp}(\mathfrak{m}_{s}({\mathbb F}_{q}))$ et $f'\in C_{cusp}(\mathfrak{m}'_{s'}({\mathbb F}_{q}))$. A la suite de bien d'autres auteurs, on a d\'ecrit dans \cite{W} 3.7 une fa\c{c}on non canonique de "relever" $s$ en une facette ${\cal F}(s)$ de $Imm(G_{AD})$. On ne rappelle pas la d\'efinition. La propri\'et\'e essentielle de ${\cal F}(s)$ est que l'espace $\mathfrak{g}_{{\cal F}(s)}$ s'identifie naturellement \`a $\mathfrak{m}_{s}$. En particulier, $f$ s'identifie \`a un \'el\'ement de $C_{cusp}(\mathfrak{g}_{{\cal F}(s)}({\mathbb F}_{q}))$, ce qui permet de d\'efinir la fonction $f_{{\cal F}(s)}\in C_{c}^{\infty}(\mathfrak{g}(F))$. Notons $N(M,s,M',s')$ l'ensemble des $g\in G(F)$ tels que $gMg^{-1}=M'$ et $gs=s'$. Cet ensemble peut \^etre vide. Il est invariant \`a gauche par $A_{M'}(F)K_{s'}^0$, on note $\underline{N}(M,s,M',s')$ un ensemble de repr\'esentants du quotient $A_{M'}(F)K_{s'}^0\backslash N(M,s,M',s')$. Enfin, on note  $A_{G}(F)_{c}$ le plus grand sous-groupe compact de $A_{G}(F)$. 
 
 \begin{prop}{Sous ces hypoth\`eses, on a:
 
 (i) si $N(M,s,M',s')=\emptyset$, $D^G_{f'}(f_{{\cal F}(s)})=0$;
 
 (ii) si $N(M,s,M',s')\not=\emptyset$, 
 $$D^G_{f'}(f_{{\cal F}(s)})=mes(K^0_{{\cal F}(s)})mes(K^0_{s'})mes(A_{M'}(F)_{c})^{-1}\sum_{g\in \underline{N}(M,s,M',s')}(f',{^gf}).$$}\end{prop}
  
  On a d\'emontr\'e en \cite{W} 3.6 une proposition analogue pour des fonctions sur le groupe $G(F)$. La preuve est similaire pour les fonctions sur l'alg\`ebre de Lie.

  \section{L'espace $I_{\star}(\mathfrak{g}(F))$}\label{star}
  On note $I_{\star}(\mathfrak{g}(F))$ le sous-espace des $f\in I(\mathfrak{g}(F))$ tels que $I^G(X,f)=0$ pour $X\not\in \mathfrak{g}_{tn}(F)$ et $I^G(X,\hat{f})=0$ pour $X\not\in \mathfrak{g}_{ent}(F)$. Notons $Fac(G)$ l'ensemble des facettes de $Imm(G_{AD})$. Posons
$${\cal E}_{nil}(\mathfrak{g}(F))=\sum_{{\cal F}\in Fac(G)}C_{nil}(\mathfrak{g}_{{\cal F}}({\mathbb F}_{q})).$$
 Comme dans le paragraphe pr\'ec\'edent, les espaces intervenant sont identifi\'es \`a des sous-espaces de $C_{c}^{\infty}(\mathfrak{g}(F))$. Avec les notations de ce paragraphe, posons aussi
$${\cal E}_{nil,\star}(\mathfrak{g}(F))=\sum_{M\in \underline{{\cal L}}}\sum_{s\in \underline{S}^G(M)}C_{nil,cusp}(\mathfrak{m}_{s}({\mathbb F}_{q}))^{K_{s}^{\dag,G}}.$$
On envoie cet espace dans $C_{c}^{\infty}(\mathfrak{g}(F))$ de la fa\c{c}on (non canonique) suivante: 
pour tout Levi $M\in \underline{{\cal L}}$ et tout $s\in S^G(M)$, on fixe  une facette ${\cal F}(s)$ comme dans le paragraphe pr\'ec\'edent; alors  $C_{nil,cusp}(\mathfrak{m}_{s}({\mathbb F}_{q}))^{K_{s}^{\dag,G}}$ s'identifie \`a un sous-espace de $C_{nil,cusp}(\mathfrak{g}_{{\cal F}(s)}({\mathbb F}_{q}))$ et on envoie toute fonction $f\in C_{nil,cusp}(\mathfrak{m}_{s}({\mathbb F}_{q}))^{K_{s}^{\dag,G}}$ sur la fonction $f_{{\cal F}(s)}$.  Pour tout sous-espace $E\subset C_{c}^{\infty}(\mathfrak{g}(F))$, notons $IE$ son image naturelle dans $I(\mathfrak{g}(F))$.

  \begin{lem}{On a les \'egalit\'es $I{\cal E}_{nil,\star}(\mathfrak{g}(F))= I{\cal E}_{nil}(\mathfrak{g}(F))=I_{\star}(\mathfrak{g}(F))$. L'homomorphisme naturel ${\cal E}_{nil,\star}(\mathfrak{g}(F))\to I_{\star}(\mathfrak{g}(F))$ est bijectif.} \end{lem} 
  
  Preuve. Par d\'efinition, on a l'inclusion ${\cal E}_{nil,\star}(\mathfrak{g}(F))\subset  {\cal E}_{nil}(\mathfrak{g}(F))$. La m\^eme preuve qu'en \cite{W} lemme 3.7 (qui concernait les groupes) montre que cette inclusion devient une \'egalit\'e dans l'espace $I(\mathfrak{g}(F))$, c'est-\`a-dire $I{\cal E}_{nil,\star}(\mathfrak{g}(F))= I{\cal E}_{nil}(\mathfrak{g}(F))$. Pour ${\cal F}\in Fac(G)$ et $f\in C_{nil}(\mathfrak{g}_{{\cal F}}({\mathbb F}_{q}))$, la fonction $f_{{\cal F}}$ est \`a support topologiquement nilpotent et sa transform\'ee de Fourier $\hat{f}$ est \`a support entier. D'o\`u l'inclusion $I{\cal E}_{nil}(\mathfrak{g}(F))\subset I_{\star}(\mathfrak{g}(F)$.
  
  Pour plus de clart\'e, notons $\phi$ la transformation de Fourier dans $C_{c}^{\infty}(\mathfrak{g}(F))$ ou $I(\mathfrak{g}(F))$. L'espace $\phi({\cal H}_{0})$ est la somme sur ${\cal F}\in Fac(G)$ des sous-espaces des fonctions \`a support dans $\mathfrak{k}_{{\cal F},0+}$. Puisque la r\'eunion de ces ensembles est $\mathfrak{g}_{tn}(F)$, $\phi({\cal H}_{0})$ est l'espace des fonctions \`a support dans $\mathfrak{g}_{tn}(F)$. Soit $f\in I_{\star}(\mathfrak{g}(F)$. Puisque les int\'egrales orbitales de $f$ sont nulles hors de $\mathfrak{g}_{tn}(F)$, on a $f\in \phi(I({\cal H}_{0}))$, ou encore $\phi(f)\in I({\cal H}_{0})$. Supposons que $\phi(f)$ soit annul\'ee par $I(\mathfrak{g}(F))^*_{nil}$. La proposition \ref{debacker} dit que $\phi(f)$ est aussi annul\'ee par $I(\mathfrak{g}(F))^*_{0}$, autrement dit que les int\'egrales orbitales $I^G(X,\phi(f))$ sont nulles pour $X\in \mathfrak{g}_{ent}(F)$. Par hypoth\`ese elles sont nulles aussi pour $X\not\in \mathfrak{g}_{ent}(F)$. Donc $\phi(f)=0$. Il en r\'esulte que l'homomorphisme naturel
$$I(\mathfrak{g}(F))^*_{nil}\to \phi( I_{\star}(\mathfrak{g}(F)))^*$$
est surjectif. Il r\'esulte de ce que l'on a dit \`a la fin du paragraphe \ref{quelquesespaces} que $I{\cal E}_{nil}(\mathfrak{g}(F))=\phi(I{\cal H}_{0}^{\sharp})$ (les espaces ${\cal E}_{nil}(\mathfrak{g}(F))$ et $\phi({\cal H}_{0}^{\sharp})$ ne sont pas tout-\`a-fait \'egaux car les fonctions attach\'ees \`a une facette ${\cal F}$ sont suppos\'ees invariantes par conjugaison par $G_{{\cal F}}({\mathbb F}_{q})$ dans le cas de l'espace ${\cal E}_{nil}(\mathfrak{g}(F))$ mais cela ne fait pas de diff\'erence quand on envoie les espaces dans $I(\mathfrak{g}(F))$). L'homomorphisme ci-dessus se prolonge en la suite
$$I(\mathfrak{g}(F))^*_{nil}\to \phi( I_{\star}(\mathfrak{g}(F)))^*\to \phi(I{\cal E}_{nil}(\mathfrak{g}(F)))^*=I{\cal H}_{0}^{\sharp,*}.$$
La proposition \ref{debacker} dit que la compos\'ee est injective. Puisque le premier homomorphisme est surjectif, le second est injectif. Dualement et apr\`es transformation de Fourier, cela nous dit que l'inclusion $I{\cal E}_{nil}(\mathfrak{g}(F))\to 
I_{\star}(\mathfrak{g}(F))$ est surjective.  Cela d\'emontre la premi\`ere assertion de l'\'enonc\'e. 
  
 Dans le paragraphe pr\'ec\'edent, on a d\'efini un espace ${\cal D}(\mathfrak{g}(F))$. Il se construit \`a l'aide d'espace $C_{cusp}(\mathfrak{g}_{s}({\mathbb F}_{q}))$. En rempla\c{c}ant dans les d\'efinitions chacun de ces espaces par $C_{nil,cusp}(\mathfrak{g}_{s}({\mathbb F}_{q}))$, on construit de m\^eme un espace ${\cal D}_{nil}(\mathfrak{g}(F))$, qui est un sous-espace de ${\cal D}(\mathfrak{g}(F))$. La description que l'on a donn\'ee de cet espace et la proposition \ref{produitscalaire} entra\^{\i}nent que l'on peut fixer des bases $(f_{i})_{i=1,...,n}$ de ${\cal D}_{nil}(\mathfrak{g}(F))$ et $(\varphi_{i})_{i=1,...,n}$ de ${\cal E}_{nil,\star}(\mathfrak{g}(F))$ de sorte que la matrice $(D^G_{f_{i}}(\varphi_{j}))_{i,j=1,...,n}$ soit inversible. A fortiori, les images des $\varphi_{i}$ dans $I(\mathfrak{g}(F))$ sont lin\'eairement ind\'ependantes. Donc l'homomorphisme naturel ${\cal E}_{nil,\star}(\mathfrak{g}(F))\to I(\mathfrak{g}(F))$ est injectif. Alors la seconde assertion de l'\'enonc\'e se d\'eduit de la premi\`ere. $\square$

 \section{L'espace $FC(\mathfrak{g}(F))$}\label{FC}
On note $FC(\mathfrak{g}(F))$ le sous-espace des fonctions $f\in I(\mathfrak{g}(F))$ telles que $I^G(X,f)=I^G(X,\hat{f})=0$ pour tout $X\not\in \mathfrak{g}_{tn}(F)$. Autrement dit, en notant encore $\phi$ la transformation de Fourier, $FC(\mathfrak{g}(F))=I_{\star}(\mathfrak{g}(F))\cap \phi(I_{\star}(\mathfrak{g}(F))$. 

{\bf Remarque.} Comme on l'a dit au paragraphe \ref{corpsfinis}, l'espace $FC(\mathfrak{g}(F))$ d\'epend du groupe $G$ et pas seulement de son alg\`ebre de Lie.
\bigskip

Pour tout $s\in S(G)$, on d\'efinit l'espace $FC(\mathfrak{g}_{s}({\mathbb F}_{q}))\subset C_{nil,cusp}(\mathfrak{g}_{s}({\mathbb F}_{q}))$, cf. paragraphe  \ref{corpsfinis}. Posons
$$fc(\mathfrak{g}(F))=\sum_{s\in S(G)}FC(\mathfrak{g}_{s}({\mathbb F}_{q})).$$

\begin{prop}{On a l'\'egalit\'e $Ifc(\mathfrak{g}(F))=FC(\mathfrak{g}(F))$ et l'inclusion $FC(\mathfrak{g}(F))\subset I_{cusp}(\mathfrak{g}(F))$. }\end{prop}

 Preuve. L'inclusion $Ifc(\mathfrak{g}(F))\subset FC(\mathfrak{g}(F))$ est imm\'ediate. Soit $f\in FC(\mathfrak{g}(F))$. Puisque $FC(\mathfrak{g}(F))\subset I_{\star}(\mathfrak{g}(F))$, on peut d'apr\`es le lemme pr\'ec\'edent relever $f$ en un unique \'el\'ement $\varphi\in {\cal E}_{nil,\star}(\mathfrak{g}(F))$. Conform\'ement \`a la d\'efinition de cet espace, \'ecrivons
 $$\varphi=\sum_{M\in \underline{{\cal L}}}\sum_{s\in \underline{S}^G(M)}\varphi_{M,s},$$
 avec $\varphi_{M,s}\in C_{nil,cusp}(\mathfrak{m}_{s}({\mathbb F}_{q}))^{K_{s}^{\dag,G}}$. Pour tout couple $(M,s)$, d\'ecomposons $\hat{\varphi}_{M,s}$ en $h_{nil,M,s}+h^+_{M,s}$, o\`u $h_{nil,M,s}$ est \`a support nilpotent et le support de $h^+_{M,s}$ ne contient pas de nilpotents. Ces fonctions restent cuspidales et invariantes par $K_{s}^{\dag,G}$.   On a les \'egalit\'es
  $$D^G_{h^+_{M,s}}(\hat{f})=D^G_{h^+_{M,s}}(\hat{\varphi})=\sum_{M'\in \underline{{\cal L}}}\sum_{s'\in \underline{S}^G(M')}D^G_{h^+_{M,s}}((\hat{\varphi}_{M,s})_{{\cal F}(s)}).$$
  Appliquons la proposition \ref{produitscalaire}: tous les termes sont nuls sauf celui index\'e par $M'=M$ et $s'=s$. Pour celui-ci, on a
 $$D^G_{h^+_{M,s}}((\hat{\varphi}_{M,s})_{{\cal F}(s)})=c( h^+_{M,s},\hat{\varphi}_{M,s})=c(h^+_{M,s},h^+_{M,s}),$$
 o\`u $c$ est une certaine constante strictement positive. D'o\`u
 
(1) $D^G_{h^+_{M,s}}(\hat{f})=c(h^+_{M,s},h^+_{M,s})$. 

Le support de la  fonction $(h^+_{M,s})_{s}\in C_{c}^{\infty}(\mathfrak{m}(F))$  ne contient pas d'\'el\'ement topologiquement nilpotent car un \'el\'ement de $\mathfrak{k}_{s}\subset \mathfrak{m}(F)$ est topologiquement nilpotent si et seulement si sa r\'eduction dans $\mathfrak{m}_{s}$ est nilpotente. Il r\'esulte alors des d\'efinitions que, pour $\phi\in C_{c}^{\infty}(\mathfrak{m}(F))$, $D^M_{h^+_{M,s}}(\phi)=0$ si le support de $\phi$  est inclus dans  $\mathfrak{m}_{tn}(F)$. Cette propri\'et\'e se conserve par induction: pour $\phi\in C_{c}^{\infty}(\mathfrak{g}(F))$, $D^G_{h^+_{M,s}}(\phi)=0$ si le support de $\phi$ est inclus dans  $\mathfrak{g}_{tn}(F)$. Par hypoth\`ese, les int\'egrales orbitales de $\hat{f}$ sont \`a support dans $\mathfrak{g}_{tn}(F)$. Il en r\'esulte que l'on peut repr\'esenter $\hat{f}$ par un \'el\'ement de $C_{c}^{\infty}(\mathfrak{g}(F))$ qui est \`a support dans $\mathfrak{g}_{tn}(F)$ (on repr\'esente $\hat{f}$ par une fonction quelconque et on remplace celle-ci par son produit avec la fonction caract\'eristique de $\mathfrak{g}_{tn}(F)$). Donc $D^G_{h^+_{M,s}}(\hat{f})=0$. Avec (1), cela entra\^{\i}ne que $h^+_{M,s}=0$. Autrement dit, la fonction $\hat{\varphi}_{M,s}$ est \`a support nilpotent. Puisque $\varphi_{M,s}$ a la m\^eme propri\'et\'e, on a $\varphi_{M,s}\in FC(\mathfrak{m}_{s}({\mathbb F}_{q}))$. Si $M\not=G$, on a $Z(M_{s})^0\not=\{1\}$ et l'espace $FC(\mathfrak{m}_{s}({\mathbb F}_{q}))$ est nul d'apr\`es \ref{corpsfinis} (2). Mais alors, la fonction $\varphi$ appartient \`a $fc(\mathfrak{g}(F))$ et $f$ appartient \`a $Ifc(\mathfrak{g}(F))$. Cela d\'emontre la premi\`ere \'egalit\'e de l'\'enonc\'e. 

On a l'inclusion $fc(\mathfrak{g}(F))\subset {\cal S}_{cusp}(\mathfrak{g}(F))$ et on a dit au paragraphe \ref{produitscalaire} que $I{\cal S}_{cusp}(\mathfrak{g}(F))$ \'etait inclus dans $I_{cusp}(\mathfrak{g}(F))$. D'o\`u la seconde assertion de l'\'enonc\'e. $\square$

 On note simplement $\underline{S}(G)=\underline{S}^G(G)$, c'est-\`a-dire que $\underline{S}(G)$ est un ensemble de repr\'esentants des orbites de l'action de $G(F)$ dans $S(G)$. Comme au paragraphe  \ref{produitscalaire}, on voit que le sous-espace
 $$(2) \qquad \sum_{s\in \underline{S}(G)}FC(\mathfrak{g}_{s}({\mathbb F}_{q}))^{K_{s}^{\dag}}$$
 de $fc(\mathfrak{g}(F))$ s'envoie bijectivement sur $FC(\mathfrak{g}(F))$. 
 
 Notons $A_{G}^{nr}$ le plus grand sous-tore de $Z(G)$ qui soit d\'eploy\'e sur $F^{nr}$. On a
 
 (3) si $A_{G}^{nr}\not=\{1\}$, $FC(\mathfrak{g}(F))=\{0\}$.
 
 En effet,  soit $s\in S(G)$.  L'hypoth\`ese $A_{G}^{nr}\not=\{1\}$ entra\^{\i}ne que $Z(G_{s})^{0}\not=\{1\}$ donc aussi $FC(\mathfrak{g}_{s}({\mathbb F}_{q}))=\{0\}$ d'apr\`es \ref{corpsfinis}(2). 
 
 \section{Endoscopie}\label{endoscopie}
 Supposons $G$ quasi-d\'eploy\'e sur $F$. Pour $X\in \mathfrak{g}_{reg}(F)$ et $f\in C_{c}^{\infty}(\mathfrak{g}(F))$, on d\'efinit l'int\'egrale orbitale stable $S^G(X,f)=\sum_{X'}I^G(X',f)$, o\`u $X'$ parcourt un ensemble de repr\'esentants des classes de conjugaison par $G(F)$ dans la classe de conjugaison stable de $X$. On note $SI(\mathfrak{g}(F))$ le quotient de $C_{c}^{\infty}(\mathfrak{g}(F))$  par le sous-espace des $f$ telles que $S^G(X,f)=0$ pour tout $X\in \mathfrak{g}_{reg}(F)$. Les int\'egrales orbitales stables peuvent \^etre consid\'er\'ees comme des formes lin\'eaires sur ce quotient. 
 
 Notons $I^{st}_{cusp}(\mathfrak{g}(F))$ le sous-espace de $I_{cusp}(\mathfrak{g}(F))$ form\'e des $f\in I_{cusp}(\mathfrak{g}(F))$ dont les int\'egrales orbitales sont constantes sur les classes de conjugaison stable. On sait  que cet espace s'envoie bijectivement sur le sous-espace  $SI_{cusp}(\mathfrak{g}(F))$ de $SI(\mathfrak{g}(F))$  form\'e des $f\in SI(\mathfrak{g}(F))$ telles que $S^G(X,f)=0$ pour tout $X\in \mathfrak{g}_{reg}(F)-\mathfrak{g}_{ell}(F)$, cf. \cite{MW} proposition I.4.11. On identifie ces deux espaces. 
 
 Revenons au cas g\'en\'eral o\`u $G$ n'est pas suppos\'e quasi-d\'eploy\'e. On fixe un ensemble $Endo_{ell}(G)$ de classes d'\'equivalence de donn\'ees endoscopiques elliptiques de $G$. On adopte pour ces donn\'ees les notations et d\'efinitions de \cite{MW} chapitre I. Soit ${\bf G}'=(G',s,{\cal H})\in Endo_{ell}(G)$. On fixe un facteur de transfert $\Delta^{{\bf G}'}$ sur $\mathfrak{g}'(F)\times \mathfrak{g}(F)$ (pr\'ecis\'ement, il est d\'efini sur un sous-ensemble d\'etermin\'e par des conditions de r\'egularit\'e, on ne perd rien \`a consid\'erer ici qu'il est nul en dehors). Cela permet de d\'efinir l'application de transfert
 $$transfert^{{\bf G}'}:I(\mathfrak{g}(F))\to SI(\mathfrak{g}'(F)).$$
 Le groupe d'automorphismes $Aut({\bf G}')$ agit naturellement sur l'espace d'arriv\'ee et l'image du transfert est contenu dans le sous-espace des invariants.
 
 {\bf Remarque.} On utilise la d\'efinition de \cite{MW} I.2.6 de cette action du groupe d'automorphismes, qui n'est pas celle que l'on trouve dans d'autres r\'ef\'erences. Un automorphisme $x\in Aut({\bf G}')$ d\'etermine un automorphisme alg\'ebrique $\alpha_{x}$ de $G'$, uniquement d\'efini \`a automorphismes int\'erieurs pr\`es. Il existe $c\in {\mathbb C}^{\times}$
tel que, pour tous $(X',X)\in \mathfrak{g}'(F)\times \mathfrak{g}(F)$, on ait $\Delta^{{\bf G}'}(\alpha_{x}(X'),X)=c\Delta^{{\bf G}'}(X',X)$.   
  Pour $f\in SI(\mathfrak{g}'(F))$, l'image de $f$ par $x$ est alors la fonction $x(f)$ d\'efinie par   l'\'egalit\'e $x(f)(X')=cf(\alpha_{x}^{-1}(X'))$ pour tout $X'\in \mathfrak{g}'(F)$.   La constante $c$  vaut $1$ si $G$ est quasi-d\'eploy\'e.   
 \bigskip

 Notons $I_{cusp}(\mathfrak{g}(F),{\bf G}')$ le sous-espace des $f\in I_{cusp}(\mathfrak{g}(F))$ telles que $transfert^{{\bf G}''}(f)=0$ pour tout ${\bf G}''\in Endo_{ell}(G)-\{{\bf G}'\}$. Alors
   le transfert  se restreint en un isomorphisme
 de $I_{cusp}(\mathfrak{g}(F),{\bf G}')$ sur $I^{st}_{cusp}(\mathfrak{g}'(F))^{Aut({\bf G}')}$. Ce dernier est une similitude pour les produits scalaires elliptiques, cf. \cite{MW} I.4.17. Pr\'ecis\'ement, posons
 $$c(G,{\bf G}')=\vert \pi_{0}(Z(\hat{G})^{\Gamma_{F}})\vert \vert \pi_{0}(Z(\hat{G}')^{\Gamma_{F}})\vert ^{-1}\vert Aut({\bf G}')\vert ^{-1}.$$
 Soient $f_{1},f_{2}\in I_{cusp}(\mathfrak{g}(F),{\bf G}')$, notons $f'_{i}=transfert^{{\bf G}'}(f_{i})$ pour $i=1,2$. Alors
 $$(1) \qquad J^G_{ell}(f_{1},f_{2})=c(G,{\bf G}')J^{G'}_{ell}(f'_{1},f'_{2}).$$
 On a fix\'e une forme bilin\'eaire sym\'etrique non d\'eg\'en\'er\'ee $<.,.>$ sur $\mathfrak{g}(F)$. Il lui correspond une forme analogue sur $\mathfrak{g}'(F)$, que l'on note encore $<.,.>$, qui est caract\'eris\'ee par la propri\'et\'e suivante. Soit $T'$ un sous-tore maximal elliptique de $G'$. On sait qu'il existe un sous-tore maximal elliptique $T$ de $G$ et un isomorphisme $\xi:T\to T'$ d\'efini sur $F$ tel que, pour $X\in \mathfrak{t}(F)\cap \mathfrak{g}_{reg}(F)$, la classe de conjugaison stable de $X$ correspond \`a celle de $X'=\xi(X)$. On a alors $<X,Y>=<\xi(X),\xi(Y)>$ pour tous $X,Y\in \mathfrak{t}(F)$. Cette forme bilin\'eaire sur $\mathfrak{g}'(F)$ a les m\^emes propri\'et\'es que la forme initiale sur $\mathfrak{g}(F)$. A l'aide de cette forme, on d\'efinit la transformation de Fourier dans $C_{c}^{\infty}(\mathfrak{g}'(F))$. Il existe une constante $\gamma_{\psi}(G,G')\in {\mathbb C}^{\times}$ v\'erifiant $\gamma_{\psi}(G,G')^4=1$ et telle que, pour tout $f\in I(\mathfrak{g}(F))$, on ait 
 
 (2) $transfert^{{\bf G}'}(\hat{f})=\gamma_{\psi}(G,G')(transfert^{{\bf G}'}(f))\hat{}$.

On a l'\'egalit\'e 
 $$(3)\qquad I_{cusp}(\mathfrak{g}(F))=\oplus_{{\bf G}'\in Endo_{ell}(G)}I_{cusp}(\mathfrak{g}(F),{\bf G}')$$
 et cette d\'ecomposition est orthogonale pour le produit scalaire elliptique.  Il y a dans $Endo_{ell}(G)$ une donn\'ee "principale" ${\bf G}=(G^*,1,{^LG})$, o\`u $G^*$ est la forme int\'erieure quasi-d\'eploy\'ee de $G$. Si $G$ est quasi-d\'eploy\'e,  on a l'\'egalit\'e $I^{st}_{cusp}(\mathfrak{g}(F))=I_{cusp}(\mathfrak{g}(F),{\bf G})$. 
 
 \section{Endoscopie et espaces $FC(\mathfrak{g}(F))$}\label{laproposition}
 Pour tout ${\bf G}'\in Endo_{ell}(G)$, posons $FC(\mathfrak{g}(F),{\bf G}')=FC(\mathfrak{g}(F))\cap I_{cusp}(\mathfrak{g}(F),{\bf G}')$. Dans le cas o\`u $G$ est quasi-d\'eploy\'e, notons $FC^{st}(\mathfrak{g}(F))=FC(\mathfrak{g}(F),{\bf G})$.
 
 \begin{prop}{(i) On a l'\'egalit\'e
 $$FC(\mathfrak{g}(F))=\oplus_{{\bf G}'\in Endo_{ell}(G)}FC(\mathfrak{g}(F),{\bf G}').$$
 
 (ii) Supposons $G$ quasi-d\'eploy\'e. L'espace $FC^{st}(\mathfrak{g}(F))$ est invariant par tout automorphisme de $G$. L'image de $FC^{st}(\mathfrak{g}(F))$ dans $SI(\mathfrak{g}(F))$ est le sous-espace des $f\in SI(\mathfrak{g}(F))$ telles que $SI^G(X,f)=SI^G(X,\hat{f})=0$ pour tout $X\in \mathfrak{g}_{reg}(F)$ qui n'est pas topologiquement nilpotent.
 
 (iii) Pour ${\bf G}'\in Endo_{ell}(G)$, le transfert se restreint en un isomorphisme de $FC(\mathfrak{g}(F),{\bf G}')$ sur $FC^{st}(\mathfrak{g}'(F))^{Aut({\bf G}')}$.}\end{prop}
 
 Preuve. Rappelons la propri\'et\'e basique du transfert endoscopique. Soit $X\in \mathfrak{g}_{ell}(F)$. L'ensemble des classes de conjugaison par $G(F)$ contenues dans la classe de conjugaison stable de $X$ est naturellement muni d'une structure d'espace homog\`ene principal sous un certain groupe ab\'elien fini $K(X)$. On choisit une famille $(X_{k})_{k\in K(X)}$ de repr\'esentants de ces classes de sorte que $X_{0}=X$ et que cette action soit la translation sur l'ensemble d'indices. Notons $K(X)^{\vee}$ le groupe des caract\`eres de $K(X)$. Soit $f\in I(\mathfrak{g}(F))$. Pour $\kappa\in K(X)^{\vee}$, on pose
 $$J^{\kappa}(X,f)=\sum_{k\in K(X)}\kappa(k)J^G(X_{k},f).$$
 Il existe un unique ${\bf G}'_{\kappa}\in Endo_{ell}(G)$ et un \'el\'ement $X'_{\kappa}\in \mathfrak{g}'_{ell}(F)$ de sorte que, pour tout $f\in I(\mathfrak{g}(F))$, on ait
 $$(1) \qquad SI^{G'_{\kappa}}(X'_{\kappa},transfert^{{\bf G}'_{\kappa}}(f))=\Delta^{{\bf G}'}(X'_{\kappa},X)J^{\kappa}(X,f).$$
 
 {\bf Remarque.} L'\'el\'ement $X'_{\kappa}$ n'est pas unique mais, si $X'_{\kappa}$ v\'erifie la propri\'et\'e ci-dessus, un autre \'el\'ement $Y'_{\kappa}\in   \mathfrak{g}'_{ell}(F)$ la v\'erifie aussi si et seulement s'il existe un automorphisme $\alpha$ de $G'$ associ\'e \`a un \'el\'ement de $Aut({\bf G}')$ tel que $Y'_{\kappa}$ soit stablement conjugu\'e \`a $\alpha(X'_{\kappa})$. 
 
 \bigskip
 Par inversion de Fourier sur le groupe $K(X)$, on a l'\'egalit\'e
 $$J^G(X,f)=\vert K(X)\vert ^{-1}\sum_{\kappa\in K(X)^{\vee}}J^{\kappa}(X,f).$$
 Supposons $f\in I_{cusp}(\mathfrak{g}(F))$ et \'ecrivons $f=\sum_{{\bf G}'\in Endo_{ell}(G)}f_{{\bf G}'}$ conform\'ement \`a la d\'ecomposition \ref{endoscopie}(3). Pour tout ${\bf G}'\in Endo_{ell}(G)$, on a alors l'\'egalit\'e
 $$(2) \qquad J^G(X,f_{{\bf G}'})=\vert K(X)\vert ^{-1}\sum_{\kappa\in K(X)^{\vee}, {\bf G}'_{\kappa}={\bf G}'}J^{\kappa}(X,f).$$
 
  Il r\'esulte de ces formules que les int\'egrales orbitales de $f$ sont \`a support topologiquement nilpotent si et seulement si il en est de m\^eme pour toutes les fonctions $f_{{\bf G}'}$. D'autre part, d'apr\`es \ref{endoscopie}(2), on a l'\'egalit\'e $(\hat{f})_{{\bf G}'}=(f_{{\bf G}'})\hat{}$ pour tout ${\bf G}'$. Donc les int\'egrales orbitales de $\hat{f}$ sont \`a support topologiquement nilpotent si et seulement s'il en est de m\^eme pour toutes les fonctions $(f_{{\bf G}'})\hat{}$. Ces deux propri\'et\'es r\'eunies entra\^{\i}nent l'assertion (i) de l'\'enonc\'e.
  
  La premi\`ere propri\'et\'e du (ii) est imm\'ediate, les espaces $FC(\mathfrak{g}(F))$ et $I^{st}_{cusp}(\mathfrak{g}(F))$ \'etant tous deux invariants par automorphismes. Notons pour le temps de cette d\'emonstration $SI^{\natural}(\mathfrak{g}(F))$ le sous-espace des $f\in SI(\mathfrak{g}(F))$ telles que $SI^G(X,f)=SI^G(X,\hat{f})=0$ pour tout $X\in \mathfrak{g}(F)$ qui n'est pas topologiquement nilpotent. Posons $SI^{\natural}_{cusp}(\mathfrak{g}(F))=SI^{\natural}(\mathfrak{g}(F))\cap I^{st}_{cusp}(\mathfrak{g}(F))$. 
 L'inclusion $FC^{st}(\mathfrak{g}(F))\subset SI^{\natural}_{cusp}(\mathfrak{g}(F))$ est claire.   Inversement, soit $f\in SI^{\natural}_{cusp}(\mathfrak{g}(F))$. Comme on l'a dit, on consid\`ere $f$ comme un \'el\'ement de $I^{st}_{cusp}(\mathfrak{g}(F))$, a fortiori comme un \'el\'ement de $I_{cusp}(\mathfrak{g}(F))$. Pour $X\in \mathfrak{g}_{ell}(F)$, on a l'\'egalit\'e $J^G(X,f)=\vert K(X)\vert ^{-1}S^G(X,f)$. Puisque les int\'egrales orbitales stables de $f$ sont \`a support topologiquement nilpotent, ses int\'egrales orbitales le sont aussi. Il en est de m\^eme pour $\hat{f}$. Donc $f\in FC^{st}(\mathfrak{g}(F))$, ce qui d\'emontre l'\'egalit\'e 
 $$FC^{st}(\mathfrak{g}(F))= SI^{\natural}_{cusp}(\mathfrak{g}(F)).$$
 Il nous reste \`a prouver l'\'egalit\'e 
 
 (3) $SI^{\natural}(\mathfrak{g}(F))=SI^{\natural}_{cusp}(\mathfrak{g}(F))$. 
 
 Il y a une filtration $({\cal F}^nSI(\mathfrak{g}(F)))_{n=n_{0},..,n_{1}}$ de l'espace $SI(\mathfrak{g}(F))$, cf. \cite{MW} I.4.15.  On a $n_{0}=dim(A_{G})-1$ et $n_{1}$ est la dimension d'un plus grand sous-tore de $G$ d\'eploy\'e sur $F$. On a ${\cal F}^{n_{0}}SI(\mathfrak{g}(F))=\{0\}$ et ${\cal F}^{n_{0}+1}SI(\mathfrak{g}(F))=I^{st}_{cusp}(\mathfrak{g}(F))$. 
 En g\'en\'eral, pour $n=n_{0},...,n_{1}$, ${\cal F}^nSI(\mathfrak{g}(F))$ est le sous-espace des $f\in SI(\mathfrak{g}(F))$ tels que $SI^G(X,f)=0$ pour tout $X\in\mathfrak{g}_{reg}(F)$ tel que $dim(A_{Z_{G}(X)})>n$. Le quotient ${\cal F}^nSI(\mathfrak{g}(F))/{\cal F}^{n-1}SI(\mathfrak{g}(F))$ s'envoie injectivement dans l'espace
 $$\sum_{M}I^{st}_{cusp}(\mathfrak{m}(F)),$$
 o\`u $M$ parcourt les Levi de $G$ tels que $dim(A_{M})=n$ (pr\'ecis\'ement, le groupe $G(F)$ agit naturellement sur cet espace et l'image de l'application ci-dessus est exactement le sous-espace des invariants). Posons ${\cal F}^nSI^{\natural}(\mathfrak{g}(F))= {\cal F}^nSI(\mathfrak{g}(F))\cap SI^{\natural} (\mathfrak{g}(F))$. Il r\'esulte de la d\'efinition de l'application ci-dessus, cf \cite{MW}1.4.15, que  le quotient ${\cal F}^nSI^{\natural}(\mathfrak{g}(F))/{\cal F}^{n-1}SI^{\natural}(\mathfrak{g}(F))$ s'envoie injectivement dans l'espace
 $$\sum_{M}SI^{\natural}_{cusp}(\mathfrak{m}(F)),$$
 qui n'est autre que
 $$\sum_{M}FC^{st}(\mathfrak{m}(F))$$
 d'apr\`es ce que l'on a d\'ej\`a d\'emontr\'e. Mais, si $n>0$, a fortiori si $n> dim(A_{G})$, ces espaces sont nuls d'apr\`es \ref{FC}(3) . Donc ${\cal F}^nSI^{\natural}(\mathfrak{g}(F))={\cal F}^{n-1}SI^{\natural}(\mathfrak{g}(F))$. Par r\'ecurrence sur $n$, cela entra\^{\i}ne l'\'egalit\'e (3) et ach\`eve la preuve de (ii).

Soient ${\bf G}'\in Endo_{ell}(G)$, $f\in I_{cusp}(\mathfrak{g}(F),{\bf G}')$ et $f'=transfert^{{\bf G}'}(f)$. Une int\'egrale orbitale stable $S^{G'}(X',f')$ est par d\'efinition combinaison lin\'eaire d'int\'egrales orbitales $I^G(X,f)$ pour des $X\in \mathfrak{g}(F)$ correspondant \`a $X'$. Inversement, les formules (1) et (2) montrent qu'une int\'egrale orbitale $I^G(X,f)$ est combinaison lin\'eaire d'int\'egrales orbitales stables $S^{G'}(X',f')$ pour des $X'\in \mathfrak{g}'(F)$ correspondant \`a $X$. Pour $X$ et $X'$ se correspondant, il est imm\'ediat que $X$ est  topologiquement nilpotent si et seulement s'il en est de m\^eme de $X'$. Donc les int\'egrales orbitales de $f$ sont \`a support topologiquement nilpotent si et seulement si les int\'egrales orbitales stables de $f'$ le sont. Il en est de m\^eme pour les transforme\'ees de Fourier d'apr\`es  \ref{endoscopie}(1).  Cela entra\^{\i}ne l'assertion (iii) de l'\'enonc\'e. $\square$

  {\bf Remarque.} On a remarqu\'e en \ref{corpsfinis} que des  notations telles que $I(\mathfrak{g}(F))$ ou $FC(\mathfrak{g}(F))$ manquaient de pr\'ecision car ces objets d\'ependent du groupe $G$ et pas seulement de son alg\`ebre de Lie. 
  Mais, dans le cas o\`u  $G$ est quasi-d\'eploy\'e, les objets stables $I^{st}(\mathfrak{g}(F))$ ou $FC^{st}(\mathfrak{g}(F))$ ne d\'ependent bien que de cette alg\`ebre. Cela r\'esulte du fait qu'une distribution stable est invariante par l'action du groupe adjoint. Enon\c{c}ons une cons\'equence, o\`u nous ajoutons 
   des exposants $G$ pour pr\'eciser quel est le groupe concern\'e. On suppose $G$ quasi-d\'eploy\'e. On a d\'ej\`a dit que, si $A_{G}^{nr}\not=\{1\}$, on avait $FC(\mathfrak{g}(F))=\{0\}$, a fortiori $FC^{st}(\mathfrak{g}(F))=\{0\}$. Supposons $A_{G}^{nr}=\{1\}$. Alors de la d\'ecomposition $\mathfrak{g}=\mathfrak{z}(G)\oplus \mathfrak{g}_{AD}$ se d\'eduit un isomorphisme 
   $$FC^{G,st}(\mathfrak{g}(F))=FC^{Z(G)^0}(\mathfrak{z}(G)(F))\otimes_{{\mathbb C}} FC^{G_{AD},st}(\mathfrak{g}_{AD}(F))$$
   et l'espace 
   $FC^{Z(G)^0}(\mathfrak{z}(G)(F))$ est la droite port\'ee par la fonction caract\'eristique de l'ensemble $\mathfrak{z}(G)_{tn}(F)$. On peut \'evidemment remplacer ci-dessus $G_{AD}$ par $G_{SC}$.
   
   \section{Traduction en termes de distributions}\label{traduction}
   Dans ce paragraphe, on renforce la condition reliant $p$ et $G$ en imposant la condition $(Hyp)_{endo}(G)$ d\'efinie en \cite{W} 7.1. Disons seulement ici que cette condition est de la forme
   $$(1) \qquad p\geq c(G)(2+val_{F}(p))$$
   o\`u $c(G)$ est une certaine constante d\'ependant de $G$. 
   
   Supposons $G$ quasi-d\'eploy\'e. L'espace dual $SI(\mathfrak{g}(F))^*$ s'identifie au sous-espace des $D\in I(\mathfrak{g}(F))^*$ telles que $D(f)=0$ pour toute fonction $f\in I(\mathfrak{g}(F))$ dont la projection naturelle dans $SI(\mathfrak{g}(F))$ est nulle. Revenons au cas g\'en\'eral.
   Pour une donn\'ee endoscopique ${\bf G}'\in Endo_{ell}(G)$, l'application de transfert $transfert^{{\bf G}'}:I(\mathfrak{g}(F))\to SI(\mathfrak{g}'(F))$ se dualise en une application de transfert $transfert^{{\bf G}',*}:SI(\mathfrak{g}'(F))^*\to I(\mathfrak{g}(F))^*$.
   
   \begin{prop}{(i) Supposons $G$ quasi-d\'eploy\'e, soit $f\in FC^{st}(\mathfrak{g}(F))$. Alors $D^G_{f}\in SI(\mathfrak{g}(F))^*$.
   
   (ii) Soit ${\bf G}'\in Endo_{ell}(G)$ et $f\in FC(\mathfrak{g}(F),{\bf G}')$. Posons $f'=transfert^{{\bf G}'}(f)\in FC^{st}(\mathfrak{g}'(F))^{Aut({\bf G}')}$. Alors $c(G,{\bf G}')transfert^{{\bf G}',*}(D^{G'}_{f'})= D^G_{f}$.}\end{prop}
   
   Preuve. Sous l'hypoth\`ese (1), l'exponentielle est un isomorphisme de $\mathfrak{g}_{tn}(F)$ sur l'ensemble $G_{tu}(F)$ des \'el\'ements topologiquement unipotents de $G(F)$. On peut relever nos fonctions et distributions sur $\mathfrak{g}(F)$ en fonctions et distributions d\'efinies sur $G_{tu}(F)$ et appliquer \`a celles-ci les r\'esultats d\'emontr\'es en \cite{W} 8.3 et 8.5. Il en r\'esulte que les assertions de la proposition sont v\'erifi\'ees si et seulement si leurs "restrictions \`a $\mathfrak{g}_{ell}(F)$" le sont. 
   
   Pr\'ecis\'ement, pour l'assertion (i), il suffit de prouver que $D_{f}^G(\varphi)=0$ pour tout $\varphi\in I_{cusp}(\mathfrak{g}(F))$ dont l'image dans $SI(\mathfrak{g}(F))$ est nulle. D'apr\`es la preuve de la formule des traces locale d'Arthur, on a l'\'egalit\'e $D^G_{f}(\varphi)=J^G_{ell}(f,\varphi)$ (cela r\'esulte de l'appartenance \`a $I_{cusp}(\mathfrak{g}(F))$ des deux fonctions $f$ et $\varphi$). Les int\'egrales orbitales elliptiques de $f$ sont constantes sur les classes de conjugaison stable tandis que les int\'egrales orbitales stables de $\varphi$ sont nulles. Il  r\'esulte alors directement de la d\'efinition  \ref{groupespadiques} (2) que $J^G_{ell}(f,\varphi)=0$, ce qui d\'emontre (i).

   Pour l'assertion (ii), il suffit de prouver que $c(G,{\bf G}')transfert^{{\bf G}',*}(D^{G'}_{f'})(\varphi)=D^G_{f}(\varphi)$ pour tout $\varphi\in I_{cusp}(\mathfrak{g}(F))$. Cela \'equivaut \`a $c(G,{\bf G}')D^{G'}_{f'}(\varphi')=D^G_{f}(\varphi)$, o\`u $\varphi'=transfert^{{\bf G}'}(\varphi)$.  
   Pour la m\^eme raison que ci-dessus, cela \'equivaut \`a $c(G,{\bf G}')J^{G'}_{ell}(f',\varphi')=J^G_{ell}(f,\varphi)$. On peut supposer $\varphi\in I_{cusp}(\mathfrak{g}(F),{\bf G}'')$ pour une certaine donn\'ee ${\bf G}''\in Endo_{ell}(G)$. Si ${\bf G}''\not={\bf G}'$, on a $\varphi'=0$ et $D^G_{f}(\varphi)=0$ puisque $f\in  FC(\mathfrak{g}(F),{\bf G}')$ et que cet espace est orthogonal \`a  $I_{cusp}(\mathfrak{g}(F),{\bf G}'')$. Alors $J^{G'}_{ell}(f',\varphi')=0=J^G_{ell}(f,\varphi)$. Si ${\bf G}''={\bf G}'$, l'\'egalit\'e \`a d\'emontrer n'est autre que \ref{endoscopie} (2). Cela ach\`eve la preuve. $\square$

\end{document}